\documentclass[preprint,11pt,fleqn,times]{elsarticle}
\usepackage{amsmath}
\usepackage{amssymb}
\usepackage{theorem}
\usepackage{xspace}
\usepackage{epsfig}

\def\pth#1{\left(#1\right)}
\def\acc#1{\left\{#1\right\}}
\def\cro#1{\left[#1\right]}

\newcommand{\R}{\mathbb{R}}

\def\ebo{\textrm{\mathversion{bold}$\mathbf{\beta^0}$\mathversion{normal}}}
\def\eb{\textrm{\mathversion{bold}$\mathbf{\beta}$\mathversion{normal}}}  
\def\eth{\textrm{\mathversion{bold}$\mathbf{\theta}$\mathversion{normal}}}

\def\eR{I\!\!R}
\def\eE{I\!\!E}
\def\eP{I\!\!P}
\def\e1{1\!\!1}

\def\ef{{\bf \overset{.}{f}}}
\def\eff{{\bf \overset{..}{f}}}

\def\x{{\bf {x}}}
\def\X{{\bf {X}}}
\def\B{{\bf {B}}}
\def\A{{\bf {A}}}

\def\c{{\bf {c}}}
\theoremstyle{plain}

\newtheorem{theorem}{Theorem}[section]
\newtheorem{corollary}{Corollary}[section]
\newtheorem{lemma}{Lemma}[section]
\newtheorem{remark}{Remark}
\newtheorem{proposition}{Proposition}[section]
\newcommand{\beqn}{\begin{eqnarray*}}
\newcommand{\eeqn}{\end{eqnarray*}}

\def\id{\hbox{{\rm\bf 1}\kern-.4em\hbox{{\rm\bf 1}}}\! }

\newcommand{\N}{\mathbb{N}}

\def\argmin{\mathop{\mathrm{arg\,min}}} 
\begin{document}
%--------------------------------------------------
\begin{frontmatter}

\title{Two tests for sequential detection of a change-point in a  nonlinear model}
\author{Gabriela CIUPERCA \footnote{{\it email: Gabriela.Ciuperca@univ-lyon1.fr}} }
\address{Universit\'e de Lyon, Universit\'e Lyon 1, 
CNRS, UMR 5208, Institut Camille Jordan, 
Bat.  Braconnier, 43, blvd du 11 novembre 1918, 
F - 69622 Villeurbanne Cedex, France}

\begin{abstract}
In this paper, two  tests, based on CUSUM of the residuals  and least squares estimation, are studied to detect in real time a change-point in a nonlinear model. A first test statistic is proposed by extension of a method already used in the literature but for the linear models. It is tested the null hypothesis, at each sequential observation,  that there is no change in the model against a change presence. The asymptotic distribution of the test statistic under the null hypothesis is given and its convergence in probability to infinity is proved when a change occurs. These results will allow to build an asymptotic critical region.  Next, in order to decrease the type I error probability, a bootstrapped critical value is proposed  and a modified test is studied in a similar way.\\ 
Simulation results, using Monte-Carlo technique, for nonlinear models which have numerous applications, investigate the properties of the two statistic tests.  \\

\end{abstract}
\end{frontmatter}

\section{Introduction}
Our aim is the construction of  a test for detecting a change in a  parametric nonlinear model $Y_i=f(X_i;\eb_i)+\varepsilon_i$, $i=1, \cdots, n$. The parameter $\eb$  will be first  estimated by a parametric method and  hypothesis test  will be afterwards made by two nonparametric  statistics. The test statistics  we are going to consider are based on sequential empirical processes of parametrically estimated residuals. This problem appears in various fields, especially biology (for example: growth model or compartmental model), chemistry, industry (quality control), finance, ...   \\
Generally, there are two types of change-point problem: \textit{a posteriori} and \textit{a priori}(sequential). The \textit{a posteriori} change-point problem arises when the data are completely known at the end of the experiment to process. For this model we begins by finding the change-points number; after that their locations and the regression parameters on each interval are estimated. 
In the case of a parametric a posteriori model with change-points we can give the following references: for a constant model with $K$ change-points, a consistent estimator for $K$ was proposed by Yao and Au (1988), using the least squares  estimation method. If the errors are strongly mixing or long-range-dependent processes, always for a constant model, Lavielle and Moulines (2000) estimate the change-point number using a penalized least-squares approach. 
Bai (1999)  proposes  a test based on the likelihood for a linear model. Again, concerning the detection of a change  in a linear model we can remind papers based on information criterion of Osorio and Galea (2005), Wu (2008)  or still Nosek (2010). In a linear model, but with long memory errors, Belkhouja and Boutahar (2009) use several methods to detect the break number: three information criteria, a sequential parametric test and a procedure based on sum of squared residuals. A large class of time series with change-points are estimated by a semi-parametric framework, but for a known change number, by Bardet et al. (2012). 
 For a parametric nonlinear model, with multiple change-points, a general criterion is proposed by Ciuperca (2011).  For the detection of the change-point number by hypothesis test in a linear a posteriori model, we can remind the paper of Liu et al. (2008), where the empirical likelihood test was considered in the particular case to detect a single change in a linear model.  Qu and Perron (2007) propose likelihood ratio type statistics to test the null hypothesis $K$ changes,  against the alternative hypothesis of $K+1$ changes, always for a linear model.\\
 In the sequential change-point problem, which will be presented here, the detection is performed in real time.  In a  linear model, the  most used technique is the CUSUM method. Horv\'ath et al. (2004) propose two schemes to detect a change in a linear model, results which are improved, using the bootstrapping, by Hu\v{s}kov\'a and Kirch (2012). The same method we find in Xia et al. (2009) for a generalized linear model. \\
 In the sequential change-point detection literature most researches consider the detection of a change in the random variable distribution (see e.g. Lai and Xing, 2010, or Mei, 2006). We can also recall several testing procedures proposed by  Neumeyer and Van Keilegom (2009) for detecting the change-points in the error distribution of non-parametric regression models. \\
 In this paper, the real time change-point detection in a nonlinear model is studied. Generalizing Horvath et al. (2004) framework, a first test statistic is studied using the weighted CUSUM method, calculated after that the model parameters have been estimated be least squares method. Next, in order to decrease the type I error probability, following the idea introduced by  Hu\v{s}kov\'a and Kirch (2012) for the linear case,  a modified test (of the first) by bootstrapping is considered. It is important to note that, the nonlinearity  changes the results and the  approach  made  by Horvath et al. (2004) and by   Hu\v{s}kov\'a and Kirch (2012) for the linear case. Above all, in a linear model, the least squares estimator of the parameters has an explicit expression, which facilitates the calculations and the results proofs. All results proofs are based on the explicit form of the estimator. In the nonlinear case, since the estimator expression is unknown and  the regression function derivatives  with respect to regression parameters  depends  on parameters and on regressors as well, imply that the  theoretical results (and  their proofs) are different. These problems are even more difficult to solve in a model where change-point occurs. Numerical algorithms will also change to calculate the critical value and test  the break presence. On the other hand, in the paper of  Hu\v{s}kov\'a and Kirch (2012), the fact that the linear model contains  intercept(see the Assumption ${\cal A}$.1(ii)), influences in a important way the results.   It is worth mentioning that we don't impose a discontinuity condition in the change-point for the model.  By simulations, for two nonlinear models which have numerous practical applications, we obtain that the two proposed tests have the empirical power equal to 1 and the empirical sizes widely smaller than the fixed theoretical size. However, the precision of the change-point estimator is the same by both methods.   \\
 The paper is organized as follows. In Section 2, we introduce the model assumptions and some general notations. The construction of a statistical test and its asymptotic behavior are presented in Section 3. To decrease the type I error probability,  Section 4 presents a modified test by bootstrapping. Next, simulation results illustrate the obtained theoretical results in Section 5. The proofs of the main results are given in Section 6, followed in Appendix by some Lemmas.    
\section{Model and notations}
  For coherence,  we try to use the some notations as in Hu\v{s}kov\'a and Kirch's paper, where the linear model was considered. \\
Let us consider the following random parametric nonlinear model with independent observations
\[
Y_i=f(\X_i;\eb_{i})+\varepsilon_i, \qquad i=1,\cdots, m, \cdots, m+T_m.
\]
 For the observation $i$, $Y_i$ denotes the response variable,  $\X_i$ is a $p \times 1$ random vector of regressors, the function $f:\eR^p \times \Theta \rightarrow \eR$ is known up to the parameters  $\eb_i$ of dimension $q \times 1$, $\eb_i \in \Theta \subseteq \eR^q$, with $\Theta$ a compact set. For the function $f$ we make the classical suppositions for a nonlinear model: $\ef({\bf x}; \eb)$ is continuous in  ${\bf x}$ and of class $C^2(\Theta)$. For the function $f({\bf x};\eb)$, we denote $\ef({\bf x}; \eb) \equiv \partial f({\bf x}; \eb)/ \partial \eb$ and $\eff({\bf x};\eb) \equiv \partial^2f({\bf x};\eb)/\partial \eb^2$.  We suppose that on the first $m$ observations, no change in the parameter regression has occurred
\[
\eb_{i}=\ebo, \qquad \textrm{for } i=1,\cdots, m,
\]
with $\ebo$ the true value of the parameter on the observations $1, \cdots, m$. The value of $\ebo$ is unknown.\\
We test the null hypothesis, that for all the following observations, there is no change in the model
\begin{equation}
\label{eq3}
H_0: \eb_{i}=\ebo, \qquad m+1 \leq i \leq m+T_m,
\end{equation}
against the hypothesis that there is a change to the $m+k^0_m+1$ observation
\begin{equation}
 \label{eq4}
H_1: \exists k^0_m \geq 1,  \textrm{such that  }
 \left\{
 \begin{array}{lll}
 \eb_{i,m}=\ebo & \textrm{for } & m+1 \leq i \leq m+k^0_m \\
 \eb_{i,m}=\eb^0_m \neq \ebo & \textrm{for } & m+k^0_m+1 \leq  i \leq m+T_m.
 \end{array}
 \right.
 \end{equation}
 The value of $\eb^0_m$ is also unknown. 
 This problem has been addressed in the literature if function $f$ is linear $f({\bf x};\eb)={\bf x}^t \eb$ (see Horv\'ath et al., 2004, Hu\v{skov\'a and Kirch, 2012}).
 Let be the sequential detector statistic, built as the weighted cumulative sum of the  residuals, for $0 \leq \gamma <1/2$, $k=1, \cdots, T_m$
 \begin{equation}
 \label{eq5}
 \left\{
  \begin{array}{l}
 \Gamma(m,k,\gamma)\equiv \sum_{m+1 \leq i \leq m+k}\hat \varepsilon_i/g(m,k,\gamma)=\sum_{m+1 \leq  i \leq m+k} [ Y_i - f(\X_i,\hat \eb_m)]/g(m,k,\gamma)\\
 \textrm{with } g(m,k,\gamma) \equiv m^{1/2}\pth{1+\frac{k}{m}} \pth{\frac{k}{k+m}}^\gamma,
 \end{array}
 \right.
 \end{equation}
 where
 $ \hat \eb_m\equiv \argmin_\beta \sum^m_{j=1} [Y_j-f(\X_j;\eb)]^2 $
 is the least squares(LS) estimator of $\eb$ calculated on the observations $1, \cdots, m$. With this estimator we calculate the   parametric  residuals
$ \hat \varepsilon_i \equiv Y_i-f(\X_i;\hat \eb_m)$, for $i=1, \cdots,k$. Recall that the cumulative sum (CUSUM) of the residuals is $\sum^{m+k}_{i=m+1} \hat \varepsilon_i$.
  Let be the $ q \times q$-matrix $ \B_m\equiv m^{-1} \sum^m_{i=1} \ef(\X_i;\ebo) \ef^t(\X_i;\ebo)$ which is supposed non-regular for all $m$ with probability one.
 Classic asymptotic results for a nonlinear regression (see also the relation (\ref{estbeta})) imply
$\hat \eb_m-\ebo=\B_m^{-1} \cro{m^{-1} \sum^m_{i=1} \ef(\X_i;\ebo) \varepsilon_i}(1+o_{\eP}(1))$.
The function $g(m,k,\gamma)$ of the relation (\ref{eq5}), proposed by Horv\'ath et al.(2004), is used as a boundary.\\
 Let us also consider the notations: $\A \equiv \eE[ \ef(\X;\ebo)]$,  $\B\equiv \eE[\ef(\X;\ebo) \ef^t(\X;\ebo)]$, $A_i\equiv \ef^t(\X_i;\ebo) \B^{-1} \A$, ${\cal D}\equiv\cro{\A^t \B^{-1}\A }^{1/2}$,  $D_A \equiv \A^t \A$. Matrix $\B$ is supposed positive definite.  All throughout the paper, vectors and matrices are written in bold face.\\

 The regression function, the random vector $\X_i$ and the error $\varepsilon_i$ satisfy the following assumptions:\\
 \textbf{(A1)} $(\varepsilon_i)_{1 \leq i \leq n}$ are i.i.d. and  $\eE[\varepsilon_i]=0$, $Var[\varepsilon_i]=\sigma^2$ and $\eE[|\varepsilon_i|^\nu]< \infty$ for some $\nu >2$.\\
 \textbf{(A2)} $\eff(\x,\eb)$ is bounded for all $\eb$ in a neighborhood of $\ebo$, for all  $\x \in \eR^p$.\\
 \textbf{(A3)} For every $i=1, \cdots, T_m$, the errors  $\varepsilon_i$ are independent of the random vectors $\X_j$, for all $j=1, \cdots, m+T_m$.\\
 \textbf{(A4)} $(m+l)^{-1} \sum^{m+l}_{i=1} f(\X_i;\ebo) \overset{{a.s.}} {\underset{m \rightarrow \infty}{\longrightarrow}} \eE[f(\X;\ebo)]  $, $(m+l)^{-1} \sum^{m+l}_{i=1} \ef(\X_i;\ebo)\overset{{a.s.}} {\underset{m \rightarrow \infty}{\longrightarrow}} \eE[\ef(\X;\ebo)]  $, \\$(m+l)^{-1} \sum^{m+l}_{i=1} \ef(\X_i;\ebo) \ef^t(\X_i;\ebo)\overset{{a.s.}} {\underset{m \rightarrow \infty}{\longrightarrow}} \B  $ for all $l=0,1, \cdots, T_m$.\\
 
 Assumptions (A2) and (A4) are made for the true parameter $\ebo$, under null hypothesis $H_0$. For the parameter $\eb^0_m$, under the alternative hypothesis, we request only the similar of (A4):\\
 \textbf{(A5)} $(m+k^0_m+l)^{-1} \sum^{m+k^0_m+l}_{i=1} f(\X_i;\eb^0_m)\overset{{a.s.}} {\underset{m \rightarrow \infty}{\longrightarrow}} \eE[f(\X;\eb^0_m)]  $, $(m+k^0_m+l)^{-1} \sum^{m+k^0_m+l}_{i=1} \ef(\X_i;\eb^0_m)$ $\overset{{a.s.}} {\underset{m \rightarrow \infty}{\longrightarrow}} \eE[\ef(\X;\eb^0_m)]$, $(m+k^0_m+l)^{-1} \sum^{m+k^0_m+l}_{i=1} \ef(\X_i;\eb^0_m)\ef^t(\X_i;\eb^0_m)\overset{{a.s.}} {\underset{m \rightarrow \infty}{\longrightarrow}} \B$, for all $l=0,1, \cdots,$ $ T_m-k^0_m$.\\
 
 The assumption that the nonlinear function $f$ is continuous in ${\bf x}$, of class $C^2$ in $\eb$ and also assumptions (A2) and (A4) are commonly used in nonlinear modeling and are necessary for the consistency and the asymptotic normality of the LS parameter estimator (see e.g. Seber and Wild, 2003). Furthermore, the two values $\ebo$ and $\eb^0_m$ are interior points of the set $\Theta$.  \\ 
 The error variance $\sigma^2$ is unknown. To estimate it, on the historical observations $i=1, \cdots, m$, we consider an consistent estimator  
 \begin{equation}
 \label{eq6}
 \hat \sigma^2_m \equiv \frac{1}{m-q}\sum^m_{j=1}[Y_j-f(\X_j;\hat \eb_m)]^2. 
 \end{equation}
 For the errors, let us consider:    $\bar \varepsilon_{m+k}=(m+k)^{-1} \sum^{m+k}_{j=1}\varepsilon_j$, $ \overline{\varepsilon^2}_{m+k}=(m+k)^{-1} \sum^{m+k}_{j=1}\varepsilon^2_j$, and then,  an another estimator for its variance besides of (\ref{eq6}), built on the $m+k$ first observations, is 
$\hat \sigma^2_{m,k}\equiv (m+k)^{-1} \sum^{m+k}_{i=1}(\varepsilon_i-\bar \varepsilon_{m+k})^2$.\\
 Two cases are possible for the sample size, which will give different results, under the null hypothesis for the test statistics:
 \begin{itemize}
\item $T_m= \infty$, the open-end procedure;
\item $T_m < \infty$, $\lim_{m \rightarrow \infty} T_m=\infty$, with $\lim_{m \rightarrow \infty} \frac{T_m}{m}=T>0$, with the possibility $T=\infty$. In this case we have the closed-end procedure.
\end{itemize}
 
 Concerning the used norms, for a $p$-vector $ {\bf v}=(v_1, \cdots, v_p)$, let us denote by  $\|{\bf v}\|_1=\sum^p_{j=1}|v_j|$ its $L^1$-norm and $\|{\bf v}\|_2=(\sum^p_{j=1} v_j^2)^{1/2}$ its $L^2$-norm. For a matrix ${\bf {\cal M}}=(a_{ij})_{\overset{1\leq i \leq p}{1 \leq j \leq q}}$, we  denote by $\|{\bf {\cal M}}\|_1=\max_{j=1, \cdots, q} (\sum^p_{i=1} |a_{ij}|)$ the subordinate norm  to the vector norm  $\|. \|_1$ and by $\|{\bf {\cal M}}\|_2=\sqrt{\rho({\bf {\cal M}} {\bf {\cal M}}^t)}$  the subordinate norm to  $\|.\|_2$, with $\rho({\bf {\cal M}} {\bf {\cal M}}^t)$ the spectral radius of ${\bf {\cal M}} {\bf {\cal M}}^t$.  \\
 All throughout the paper,  $C$ denotes a positive generic constant which may take different values in different formula or even in different parts of the same formula. All vector are column and ${\bf v}^t$ denotes the transpose of ${\bf v}$. 
 We say that a random variable set $(V_n)$ is bounded by a constant $C$ with a probability close to 1 (or with a probability arbitrarily large): $\forall \epsilon >0$, $\exists n_\epsilon \in \N$ such that $\eP[ V_n >C] < 1-\epsilon$.\\

 Now, a notation and a relation on the function $g$, used for the result proofs. 
 Using the  relation that for all  $x>0$ we have $0 < \frac{x}{1+x} <1$ and that  $\gamma \in [0,1/2)$, we obtain that
 \begin{equation}
 \label{Km}
 K_m \equiv \sup_{1 \leq k < \infty} \frac{k m^{-1/2}}{g(m,k,\gamma)}=\sup_{1 \leq k < \infty} \pth{\frac{k/m}{1+k/m}}^{1-\gamma} \in [0,1].
 \end{equation} 
  
 After from these general notations, in every section we shall give the notations used for each test.\\

 The proofs of all main results of Sections 3 and 4 are given in Section 6. To prove these results, necessary lemmas  are stated and proved in Appendix (Section 7).
 
 \section{Test by weighted  CUSUM, without bootstrapping}
 
 We are going first to build  a test statistic based on the residuals $\hat \varepsilon_i =Y_i-f(\X_i; \hat \eb_m)$ after the observation $m$ by estimating the parameter $\eb$ on the historical data $(Y_i,\X_i)_{1 \leq i \leq m}$. The study of this statistic will be hampered by the fact that the estimator $\hat \eb_m$ does not have an explicit expression. \\
The following  Theorem is the generalization of the result obtained by  Horv\'ath et al.(2004) for the linear model, on the asymptotic distribution of the test statistic under the null hypothesis given by (\ref{eq3}). We remark that, unlike to the linear case, the asymptotic distribution of the test statistic, under $H_0$, depends on the function $f({\bf x};\ebo)$ and on the true parameter $\ebo$. The value of $T_m$, with respect to $m$, also influence the asymptotic distribution.

\begin{theorem}
\label{theorem 2.1}
Let us consider the assumptions  (A1)-(A4). Under the null hypothesis $H_0$ specified by (\ref{eq3}), for all real  $c>$, we have\\
(i) If $T_m=\infty$ or ($T_m <\infty$ and $\lim_{m \rightarrow \infty} T_m/m=\infty$), then
\begin{equation}
\label{PT}
\lim_{m \rightarrow \infty}\eP \cro{ \frac{1}{\hat \sigma_m} \sup_{1 \leq k <\infty}  \left|\sum^{m+k}_{i=m+1} \hat \varepsilon_i \right|/g(m,k,\gamma) \leq c}=\eP \cro{\sup_{0 \leq t \leq  \frac{1}{{\cal D}^2} }\frac{(1+t-{\cal D}^2t)|W(t)|}{t^\gamma} \leq c}.
\end{equation}
(ii) If $T_m <\infty$ and $\lim_{m \rightarrow \infty} T_m/m=T< \infty$, then the left-hand side of (\ref{PT}) is equal to
$ \eP \cro{\sup_{0 \leq t \leq  \frac{T}{1+{\cal D}^2T} }\frac{(1+t-{\cal D}^2t)|W(t)|}{t^\gamma} \leq c}$.
Here $\{ W(t), 0 \leq t < \infty\}$ is a Wiener process (Brownian motion) i.e. a centered Gaussian process, with covariance function $Cov(W(s), W(t))=min(s,t)$, $s,t \in [0,\frac{1}{{\cal D}^2}]$ for (i) and $s,t \in [0,\frac{T}{1+{\cal D}^2T}]$ for (ii).

\end{theorem}

In order to have a test statistic, thus, to build a critical region, it is necessary to study the behavior of the statistic in the left-hand side of (\ref{PT}) under the alternative hypothesis $H_1$. By the following Theorem, this statistic converges in probability to infinity as $m \rightarrow \infty$. For this, we suppose that the change-point $k^0_m$ is not very far from the last observation of historical data. Obviously, this supposition poses no problem for practical applications, since if hypothesis $H_0$ was not rejected until an observation $k_m$ of order $m$, we reconsider as historical data, all observations of 1 to $k_m$. Another supposition is that, before and after the break, on average, the model is different, without imposing a discontinuity condition in the change-point.  
\begin{theorem}
\label{theorem 2.2}
Suppose that  the assumptions  (A1)-(A5) hold.   Under the alternative hypothesis $H_1$ specified by (\ref{eq4}), if $k^0_m=O(m)$ and  $\eE[f(\X; \eb^0)] \neq \eE[f(\X;\eb^0_m)]$ hold also, then
\[
\frac{1}{\hat \sigma_m} \sup_{1 \leq k \leq T_m} \cro{\left|\sum^{m+k}_{i=m+1} \hat \varepsilon_i \right|/g(m,k,\gamma)}  \overset{{\eP}} {\underset{m \rightarrow \infty}{\longrightarrow}} \infty.
\]
\end{theorem}
Considering the Theorems \ref{theorem 2.1} and \ref{theorem 2.2} we derive in the next corollary a test statistic for testing the lack of change against the break presence. 
\begin{corollary}
\label{Corollary 1}
Consequence of these two theorems, following  statistic can be used to test   $H_0$ against $H_1$:
\begin{equation}
\label{stat_test1}
Z_\gamma(m)\equiv \frac{1}{\hat \sigma_m} \sup_{1 \leq k \leq T_m}\left|\sum^{m+k}_{i=m+1} \hat \varepsilon_i \right|/g(m,k,\gamma).
\end{equation} 
The asymptotic critical region is $\acc{ Z_\gamma(m) \geq c_\alpha(\gamma)}$, where $c_\alpha(\gamma)$ is the  $(1-\alpha)$ quantile  of the distribution of  $\sup_{0 \leq t \leq  \frac{1}{{\cal D}^2}} [ t^{-\gamma} (1+t-{\cal D}^2t)|W(t)|]$, if $\lim_{m \rightarrow \infty} T_m/m=\infty$, and of $\sup_{0 \leq t \leq  \frac{T}{1+{\cal D}^2T} }[t^{-\gamma}(1+t-{\cal D}^2t)|W(t)|]$, if $\lim_{m \rightarrow \infty} T_m/m=T \in (0,\infty)$. For some given $\alpha \in (0,1)$, this statistical test, consequence of Theorems \ref{theorem 2.1} and \ref{theorem 2.2}, has the asymptotic type I error probability (size) $\alpha$ and the asymptotic power 1. 
\end{corollary}
It is important to note that, in the linear case $f(\x;\eb)=\x^t \eb$, the value of ${\cal D}$ depends only on $\eE[\X]$, $\eE[\X \X^t]$ but not on the values of $\ebo$. For a nonlinear model, the critical values $c_\alpha(\gamma)$ depend on the regression function $f$, the distribution of random vector $\X$ and    on parameter value $\ebo$ before the change-point.
\begin{remark}
In the linear case, the assumption that the model contains intercept, $\X=(1,X_1, \cdots, X_p)$, $\eb=(b_0,b_1, \cdots, b_p)$, imposed by Horv\'ath et al. (2004), is essential. If $\eE[X_1]=\cdots = \eE[X_p]=0$, then it is necessary that the model has different intercepts before and after change-point. Without this supposition, the test statistic $Z_{\gamma}(m)$ can not converge to infinity under $H_1$.
\end{remark}

\noindent Therefore, we deduce from it that, the null hypothesis $H_0$ is rejected in the change-point
 \begin{equation}
 \label{eq7}
 \hat \tau_m \equiv \left\{
 \begin{array}{l}
 \inf \acc{1 \leq k \leq T_m, \hat \sigma_m^{-1} | \Gamma(m,k,\gamma)| \geq c_\alpha(\gamma)} \\
 \infty, \;\; \textrm{if } \hat \sigma_m^{-1} | \Gamma(m,k,\gamma)| <c_\alpha(\gamma), \;\; \textrm{for every } 1 \leq k \leq T_m.
 \end{array}
 \right.
 \end{equation}
which we can consider as estimator for $k^0_m$.

\section{Test by weighted CUSUM, with bootstrapping}
In order to improve the critical values of the test, thus, to decrease  the  type I error probability, we extend the  method proposed by  Hu\v{s}kov\'a and Kirch (2012), which uses  the bootstrapping to calculate the critical value, function of the observation position, after the observation $m$.\\ 
Let us suppose that  until the observation $m+k$, the hypothesis $H_0$ has not been rejected yet. Thus, for $l=1, \cdots, m+k$ we have that under $H_0$, using the relation (\ref{estbeta}) and  the proof of the  Lemma \ref{Lemma 5.2}, the  cumulative sum of the residuals defined by (\ref{eq5}) can be approached
\[
\Gamma(m,l,\gamma) = \cro{\sum^{m+l}_{i=m+1} \varepsilon_i- \pth{\frac{1}{m} \sum^m_{j=1} \ef^t(\X_j;\ebo) \varepsilon_j} \B^{-1}_{m}  \sum^{m+l}_{i=m+1} \ef(\X_i;\ebo)}/g(m,l,\gamma)(1+o_{\eP}(1)).
\]
In order to realize the bootstrapping, let us consider the discrete uniform random variables ${\cal U}_{m,k}(i)$, for $i=1, \cdots, m+T_m$, such that $\eP[{\cal U}_{m,k}(i)=j]=1/(m+k)$, for $j=1, \cdots, m+k$. 
We denote also  by  $\eP^*_{m,k}$, $\eE^*_{m,k}$, $Var^*_{m,k}$ the conditional probability, expectation, variance we respect to $\acc{{\cal U}_{m,k}(i), 1 \leq i \leq m+T_m }$, given $(Y_j,\X_j)_{1 \leq j \leq m+k}$. The conditional expectation with the bootstrapped  regressors is, for $i=1, \cdots, m+T_m$, 
\[
\eE^*_{m,k} [  \ef(\X_{{\cal U}_{m,k}(i)};\ebo)]=\frac{1}{m+k} \sum^{m+k}_{j=1} \ef(\X_j;\ebo).
\]
 Keeping the same notations as in the linear model of Hu\v{s}kov\'a and Kirch (2012), let us consider (see  Section 2, for the other notations), for $k=1, \cdots, T_m$, following notations
\begin{itemize}
\item $\c_1(m,k,l)\equiv  D^{-1}_A\B_{m} \left[\sum^{m+l}_{i=m+1} \ef(\X_i;\ebo) \e1_{l \leq k}+ \sum^{m+k}_{i=m+k-l+1}\ef(\X_i;\ebo) \e1_{k <l < m+ k}+l(m+k)^{-1} \sum^{m+k}_{i=1} \ef(\X_i;\ebo) \e1_{l \geq m+ k} \right]$, for $1 \leq l \leq T_m$. In the linear model, $\c_1(m,k,l)$ depends only $\X_i$.
\item $\tilde \Gamma(m,k,l,\gamma)(\varepsilon_1, \cdots, \varepsilon_{m+l})\equiv \cro{\sum^{m+l}_{i=m+1} \varepsilon_i- \pth{m^{-1} \sum^m_{j=1} \ef^t(\X_j;\ebo) \varepsilon_j} \B^{-1}_{m}  \c_1(m,k,l)}/g(m,l,\gamma)$, which is an approach of the weighted CUSUM statistic  $\Gamma(m,l,\gamma)$ given by (\ref{eq5}), in order to facilitate the bootstrap.
\item $\hat \varepsilon_{m,k}(j)\equiv Y_j-f(\X_j;\hat \eb_{m+k})$ are the residuals from the ordinary least squares method,  with $\hat \eb_{m+k} \equiv \argmin_\eb \sum^{m+k}_{j=1} [ Y_j-f(\X_j;\eb)]$.
\item $\varepsilon^*_{m,k}(i) \equiv \hat \varepsilon_{m,k}( {\cal U}_{m,k}(i))$ are the bootstrap errors.
\item $
\hat \sigma^{(*)2}_{m,k} \equiv (m-q)^{-1} \sum^{m}_{i=1} \cro{\varepsilon^*_{m,k}(i)-\pth{m^{-1}\sum^m_{j=1} \ef^t(\X_j;\ebo)\varepsilon^*_{m,k}(j)} \B^{-1}_{m} \ef(\X_i;\ebo) }^2$ the bootstrap variance estimator.\\
\item $F^*_{m,k}(x)\equiv \eP^*_{m,k} \cro{ 1/\hat \sigma^{(*)}_{m,k} \sup_{1 \leq l \leq T_m} | \tilde \Gamma(m,k,l,\gamma)( \varepsilon^*_{m,k}(1), \cdots, \varepsilon^*_{m,k}(m+l) )| \leq x }$ a distribution function calculated using the bootstrap results.
\item For $N\geq 1$, let us consider  $\tilde F_{m,k} \equiv \sum^{N-1}_{i=0} \alpha_i F^*_{m,\max((j-i)L,0)}$, for $k=jL, \cdots, (j+1)L-1$ an other distribution function, proposed by Hu\v{s}kov\'a and Kirch (2012) in order to accelerate the procedure. The positive constants $\alpha_i$ are such that $\sum^{N-1}_{i=0} \alpha_i=1$.
\end{itemize}

We note that in order to calculate the bootstrapped residuals $\varepsilon^*_{m,k}(i)$, only the data $(Y_i,\X_i)_{1 \leq i \leq m+T_m}$ are bootstrapped, not the estimator $\hat \eb_{m+k}$ of $\eb$ calculated on not bootstrapped data. \\ 
The $(1-\alpha)$ quantile  $c_{m,k;\alpha}(\gamma)$ at time $m+k$ of the distribution $\tilde F_{m,k}$ is obtained as the smallest real value such that
\begin{equation}
\label{eq19}
\tilde F_{m,k}(c_{m,k;\alpha}(\gamma)) \geq 1- \alpha.
\end{equation}
Contrary to the case of Corollary \ref{Corollary 1}, for the weighted CUSUM statistic without bootstrapping, the critical values $c_{m,k;\alpha}(\gamma)$ depend at the same time of $m$, and $k$ besides $\alpha$ and $\gamma$. \\
Before to state the main results of this section, let us recall  the H\'ajek-R\'enyi inequality (see H\'ajek and R\'enyi, 1955) that is a generalization of the Kolmogorov inequality.\\
\textit{H\'ajek-R\'enyi inequality: if ${(G_{k})}_{1 \leq k \leq n}$ is a sequence of independent random variables with $\eE[G_k]=0$, $Var(G_k)<\infty$ and  ${(b_k)}_{1 \leq k \leq n}$ is a non-decreasing sequence of positive numbers, then, for any $\epsilon>0$ and $m \leq n$, }
\[
\eP[\max_{m \leq k \leq n} \left| \frac{\sum^k_{j=1} G_j}{b_k} \right| \geq \epsilon] \leq \frac{1}{\epsilon^2} \cro{ \sum^n_{j=m+1} \frac{\eE[G^2_j]}{b^2_j} + \sum^m_{j=1} \frac{\eE[G^2_j]}{b^2_m}}.
\]
A particular case of  this inequality is we consider  $b_k=g(n,k,\gamma)$, which is an increasing sequence in  $k$, with the function $g$ specified by relation (\ref{eq5}).\\
For the linear model (see Hu\v{s}kov\'a and Kirch, 2012), to study the behavior of the distribution function $\tilde F_{m,k}$, then the behavior of the statistic $1/\hat \sigma^{(*)}_{m,k}) \sup_{1 \leq l \leq T_m} | \tilde \Gamma(m,k,l,\gamma)( \varepsilon^*_{m,k}(1), \cdots, \varepsilon^*_{m,k}(m+l) )|$, the H\'ajek-R\'enyi inequality alone was sufficient. In the nonlinear model, in the calculation of the bootstrapped residual   $\varepsilon^*_{m,k}$, then of $\hat \varepsilon_{m,k}$ , the LS estimator $\hat \eb_{m+k}$ intervenes. Since $\hat \eb_{m+k}$ was  not an explicit expression,  we need a generalization of this inequality for random variable sequence of expectation converging uniformly to 0. First, we have the following general result. 
 
\begin{proposition}
\label{Proposition A} If  ${(Z_{k,n})}_{1 \leq k \leq n}$ is a random variable such that  $\eE[Z_{k,n}]=\mu_{k,n} \rightarrow 0$, for $n \rightarrow \infty$, uniformly in $k$, and for all $\epsilon >0$ and $m \leq n$, $\eP[\max_{m \leq k \leq n} |Z_{k,n} -\mu_{k,n}| \geq \epsilon]\rightarrow 0$, then, there exists a natural number $n_\epsilon$ such that for $n \geq n_\epsilon$, $\eP[\max_{m \leq k \leq n} |Z_{k,n}| \geq 2\epsilon]\rightarrow 0$.
\end{proposition}

As a consequence of the Proposition  \ref{Proposition A} and of the H\'ajek-R\'enyi inequality, a generalization of this last one  can be established, for random variables with the expectation converging to 0. Let ${(G_j)}_{1 \leq j \leq n}$ be a sequence of random variables  such that $\eE[G^2_j] < \infty$, for all $j=1, \cdots, n$ and $\eE [ b_k^{-1} \sum^k_{j=1} G_j]=\mu_{k,n} \rightarrow 0$, uniformly in $k$, for $n \rightarrow \infty$, with the positive sequence ${(b_k)}_{1 \leq k \leq n}$ non-decreasing. Then, by the proof of Proposition \ref{Proposition A},  we have that, for any $\epsilon >0$, there exists a natural number $n_\epsilon$ such that for $n \geq n_\epsilon$ 
\begin{equation}
\label{AA}
\eP \cro{ \max_{1 \leq k \leq n} \frac{| \sum^k_{j=1} G_j|}{b_k} \geq 2 \epsilon} \leq \eP \cro{ \max_{1 \leq k \leq n} \frac{| \sum^k_{j=1} (G_j- \eE[G_j] )|}{b_k} \geq  \epsilon}.
\end{equation}
On the other hand, by the H\'ajek-R\'enyi inequality, we have for the random variable $G_j-\eE[G_j]$, for any $\epsilon >0$,
\begin{equation}
\label{BB}
\eP \cro{ \max_{1 \leq k \leq n} \frac{| \sum^k_{j=1} (G_j- \eE[G_j] )|}{b_k} \geq  \epsilon} \leq \frac{1}{\epsilon^2} \sum^n_{j=1} \frac{Var[G_j]}{b^2_j}.
\end{equation}
But $Var[G_j] \leq \eE[G^2_j]$. By the relations (\ref{AA}) and (\ref{BB}) it follows immediately that, for any sequence of random variables ${(G_j)}_{1 \leq j \leq n}$ such that $\eE[G^2_j] < \infty$, for all $j=1, \cdots, n$ and $\eE [ b_k^{-1} \sum^k_{j=1} G_j]=\mu_{k,n} \rightarrow 0$, uniformly in $k$, for $n \rightarrow \infty$, with the positive sequence ${(b_k)}_{1 \leq k \leq n}$ non-decreasing and  for any $\epsilon >0$, then, there exists a natural number $n_\epsilon$ such that for $n \geq n_\epsilon$, 
\begin{equation}
\label{CC}
\eP \cro{ \max_{1 \leq k \leq n} \frac{| \sum^k_{j=1} G_j|}{b_k} \geq 2 \epsilon} \leq \frac{1}{\epsilon^2} \sum^n_{j=1} \frac{\eE[G_j^2]}{b^2_j}.
\end{equation}

Now, in  order to study the residuals $\hat \varepsilon_{m,k}(i)=Y_i-f(\X_i;\hat \eb_{m+k})$, calculated after observation $m$, we underline, by a decomposition, the corresponding model error $\varepsilon_i$. 
Depending  on the position of the  observation  "$i$"  with respect to change-point $m+k^0_m$, where $k^0_m$ is the change-point position under the alternative hypothesis $H_1$ given by (\ref{eq4}), and on the position of $k$ with respect to $k^0_m$, we have the decomposition for the residuals
\[
\hat \varepsilon_{m,k}(i) = \varepsilon_i+f(\X_i;\ebo) \e1_{i \leq m+k^0_m} +f(\X_i;\eb^0_m) \e1_{i >m+k^0_m}- f(\X_i;\hat \eb_{m+k})  \e1_{i \leq m+k^0_m} [\e1_{k \leq k^0_m}+\e1_{k >k^0_m} 
  - f(\X_i;\hat \eb_{m+k})  \e1_{i >m+k^0_m}  \e1_{k >k^0_m} .
\]
Since $\hat \eb_{m+k}$ is the least squares estimator of $\eb$, we have 
$
0 =  \sum^{m+k}_{i=1} \ef(\X_i;\hat \eb_{m+k}) [\varepsilon_i -f(\X_i;\hat \eb_{m+k})+f(\X_i;\ebo)]\e1_{k \leq k^0_m}$  $+\sum^{m+k}_{i=m+k^0_m+1} \ef(\X_i;\hat \eb_{m+k}) [\varepsilon_i -f(\X_i;\hat \eb_{m+k})+f(\X_i;\eb^0_m)]\e1_{k > k^0_m}$. Then  $\sum^{m+k}_{i=1}\varepsilon_i \ef(\X_i;\hat \eb_{m+k}) =\sum^{m+k}_{i=1} f(\X_i;\hat \eb_{m+k}) $\\ $\cdot \ef(\X_i;\hat \eb_{m+k})
- \sum^{m+k^0_m}_{i=1} f(\X_i;\ebo) \ef(\X_i;\hat \eb_{m+k})$ 
$- \sum^{m+\min(k,k^0_m)}_{i=m+k^0_m+1} f(\X_i;\eb^0_m) \ef(\X_i;\hat \eb_{m+k})\e1_{k > k^0_m}$. 
The statistic \\ $g(m,l,\gamma) \tilde \Gamma(m,l,\gamma)( \varepsilon^*_{m,k}(1), \cdots, \varepsilon^*_{m,k}(m+l) )$ becomes
\begin{equation}
\label{uu}
\sum^{m+l}_{i=m+1} \hat \varepsilon_{m,k}({\cal U}_{m,k}(i))-   \pth{\frac{1}{m}\sum^m_{j=1} \ef^t(\X_j;\ebo)\hat \varepsilon_{m,k}({\cal U}_{m,k}(j))}  \B^{-1}_{m} \c_1(m,k,l) \equiv I_1+I_2 + {\cal R}_m, 
\end{equation}
with  $I_1 \equiv  \sum^{m+l}_{i=m+1} \varepsilon_{{\cal U}_{m,k}(i)}$ and $I_2 \equiv -   \pth{\frac{1}{m}\sum^m_{j=1} \ef^t(\X_j;\ebo) \varepsilon_{{\cal U}_{m,k}(j)} } \B^{-1}_{m} \c_1(m,k,l)$. The  expression of ${\cal R}_m$ will be specified in Appendix (Section 7). 
We precise that the bootstrapped residuals are $ \hat \varepsilon_{m,k}({\cal U}_{m,k}(i))=Y_{{\cal U}_{m,k}(i)}-f(\X_{{\cal U}_{m,k}(i)};\hat \eb_{m+k})$ and $\varepsilon_{{\cal U}_{m,k}(i)}=Y_{{\cal U}_{m,k}(i)}-f(\X_{{\cal U}_{m,k}(i)};\ebo) \e1_{{\cal U}_{m,k}(i) \leq m+k^0_m}- f(\X_{{\cal U}_{m,k}(i)};\eb^0_m) \e1_{{\cal U}_{m,k}(i) > m+k^0_m}$.
\\

With these elements, we can prove that the statistic  $\tilde \Gamma(m,k,l,\gamma)$ is  asymptotically determined  by  $I_1$ and   $I_2$ under $H_0$ and that each of them converges to a Wiener process. For these, we prove, by the following Proposition,  that the term $I_2$ can be also written asymptotically  as a sum of $\varepsilon_{{\cal U}_{m,k}(i)}$, by imposing a supplementary condition:\\
\textbf{(A6)} for any $\epsilon >0$ there exists $M>0$ such that $\eP\cro{\max_{1 \leq i \leq m} \| \ef(\X_i;\ebo) \|_2 \geq M } \leq \epsilon$.\\
The proof of Proposition \ref{Lemma 4} is given in Section 6, where the nonlinearity intervenes decisively to prove that the sum of  1 to $T_m$ for the right-hand side of an expression like (\ref{CC}) converges uniformly in probability to zero.
\begin{proposition}
\label{Lemma 4} Under the assumptions (A1)-(A4), (A6) we have for any $\epsilon >0$, in probability,
\[
\sup_{1 \leq k < \infty} \eP^*_{m,k} \cro{\max_{1 \leq l \leq T_m} \frac{\left|I_2-(-l/m \sum^m_{j=1}\varepsilon_{{\cal U}_{m,k}(j)}) \right|}{g(m,l,\gamma)} \geq \epsilon} {\underset{m \rightarrow \infty}{\longrightarrow}} 0.
\]
\end{proposition}

Taking into account the proof of Theorem \ref{theorem 2.1} concerning the asymptotic distribution of the weighted cumulative residuals sum $\Gamma(m,k, \gamma)$ calculated without bootstrapping, we show by the following results that the statistic  $\tilde \Gamma(m,k,l,\gamma)$ bootstrapped has the same asymptotic behavior under $H_0$ as $\Gamma(m,k, \gamma)$. 
Under hypothesis $H_1$, the   term ${\cal R}_m$ is asymptotically uniformly bounded and then, taking into account the relation (\ref{uu}),  $\tilde \Gamma(m,k,l,\gamma)$ is uniformly bounded a.s. also (see in Appendix,  sub-Section 7.2, the Lemmas \ref{Lemma I3I6} and \ref{Lemma I4I7I5I8}). \\

\begin{proposition}
\label{lemma 5}
Suppose that the assumptions (A1)-(A4), (A6) hold. \\
a) Under the null hypothesis $H_0$, we have, for any $x \in \R$,
\[
\sup_{1 \leq k \leq T_m} \left| \eP^*_{m,k} \cro{\frac{1}{\hat \sigma_{m,k}}\sup_{1 \leq l \leq T_m}  \tilde \Gamma(m,k,l,\gamma)(\varepsilon^*_{m,k}(1), \cdots, (\varepsilon^*_{m,k}(m+l)) \leq x}  -\eP \cro{\sup_{1 \leq l \leq T_m} \frac{|W_1 \pth{\frac{l}{m}} -\frac{l}{m}{\cal D} W_2(1) |}{\pth{1+\frac{l}{m}} \pth{\frac{l}{m+l}}^\gamma} \leq x}
\right| \overset{{\eP}} {\underset{m \rightarrow \infty}{\longrightarrow}} 0.
\] 
where $\{ W_1(t); 0 \leq t < \infty \}$ is a Wiener process, $W_2(1)$ is a standard normally distributed, independent of $\{W_1(t)\}$.\\ 
b) If furthermore the assumption (A5) holds, under the alternative  hypothesis $H_1$, for any $\epsilon >0$, there exists a constant $M>0$ such that, we have a.s.
\[
\sup_{1 \leq k \leq T_m} \eP^*_{m,k} \cro{\frac{1}{\hat \sigma_{m,k}}\sup_{1 \leq l \leq T_m}  \left| \tilde \Gamma(m,k,l,\gamma)(\varepsilon^*_{m,k}(1), \cdots, \varepsilon^*_{m,k}(m+l))\right| 
\geq M} \leq \epsilon+o_{\eP}(1).
\]
\end{proposition}

\noindent As for the Theorem  \ref{theorem 2.1}, under $H_0$, we can prove that the asymptotic distribution of $\sup_{1 \leq l \leq T_m} \frac{|W_1 (l/m) -l/m{\cal D} W_2(1) |}{(1+l/m) (l/(m+l))^\gamma}$ is $\sup_{0 \leq t \leq  \frac{1}{{\cal D}^2} }\frac{(1+t-{\cal D}^2t)|W(t)|}{t^\gamma}$ in the case $T_m=\infty$ or ($T_m <\infty$ and $\lim_{m \rightarrow \infty} T_m/m=\infty$). In the case $T_m <\infty$ and $\lim_{m \rightarrow \infty} T_m/m=T< \infty$, the asymptotic distribution is $\sup_{0 \leq t \leq  \frac{T}{1+{\cal D}^2T} }\frac{(1+t-{\cal D}^2t)|W(t)|}{t^\gamma}$. Combining Theorem \ref{theorem 2.1} with Proposition \ref{lemma 5}\textit{(a)} under the null hypothesis, on the one hand, and Theorem \ref{theorem 2.2} with Proposition \ref{lemma 5}\textit{(b)} under the alternative hypothesis, on the other hand, together with the distribution function definition $\tilde F_{m,k}$, allow to define a critical value depending of each sequential observation $k=1, \cdots, T_m$. Thus, we can define a new test statistic and study its asymptotic behavior under $H_0$ and $H_1$. 
\begin{theorem}
\label{theorem 1HK}
Suppose that the assumptions (A1)-(A4), (A6) hold and that $\alpha \in (0,1)$, $\gamma \in [0,1/2)$.  \\
a) Under the null  hypothesis $H_0$, as $m \rightarrow \infty$, we have
\[
\eP \cro{\frac{1}{\hat \sigma_m} \sup_{1 \leq k \leq T_m} \frac{|\Gamma(m,k,\gamma) |}{c_{m,k;\alpha}(\gamma)} >1 } \rightarrow \alpha.
\]
b) If furthermore the assumption (A5) holds, under the alternative hypothesis $H_1$, as $m \rightarrow \infty$, we have
\[
\eP \cro{\frac{1}{\hat \sigma_m} \sup_{1 \leq k \leq T_m} \frac{|\Gamma(m,k,\gamma) |}{c_{m,k;\alpha}(\gamma)} >1 } \rightarrow 1,
\]
with $\Gamma$ given by the relation (\ref{eq5}) and $c_{m,k;\alpha}(\gamma)$ by (\ref{eq19}) is the critical value of the distribution function $\tilde F_{m,k}$.
\end{theorem}

Thus, we are going to use as test statistic of  $H_0$, against $H_1$
\begin{equation}
\label{stat_test2}
Z_{\gamma;\alpha}^{(b)}(m) \equiv \frac{1}{\hat \sigma_m} \sup_{1 \leq k \leq T_m} \frac{|\Gamma(m,k,\gamma) |}{ c_{m,k;\alpha}(\gamma)},
\end{equation}
which will have the asymptotic critical region $\acc{Z_{\gamma;\alpha}^{(b)}(m)> 1 }$. Then the statistic $Z_{\gamma;\alpha}^{(b)}(m)$ has asymptotic size $\alpha$ and asymptotic power one for all $\gamma \in [0,1/2)$. As in Section 4, we consider the change-point estimator of $k^0_m$  is
\begin{equation}
 \label{eq7_b}
 \hat \tau_m^{(b)} \equiv \left\{
 \begin{array}{l}
 \inf \acc{1 \leq k \leq T_m,  \frac{|\Gamma(m,k,\gamma) |}{\hat \sigma_m c_{m,k;\alpha}(\gamma)} > 1}, \\
 \infty, \;\; \textrm{if }  \frac{|\Gamma(m,k,\gamma) |}{\hat \sigma_m \cdot c_{m,k;\alpha}(\gamma)} \leq 1, \;\; \textrm{for every } 1 \leq k \leq T_m.
 \end{array}
 \right.
 \end{equation}
 Then, hypothesis $H_0$ is rejected in $\hat \tau_m^{(b)} $. 
 Let us notice that, in comparison with the previous test, the value calculation $c_{m,k;\alpha}(\gamma)$ is little  more laborious, in view of the fact that, the conditional distribution functions $F^*_{m,k}$ must be first calculated.  

\section{Simulations}
In this section we report a simulation study designed  to evaluate and compare the performance of the proposed test methods. For the two methods we consider two examples: growth model and compartmental model for varied parameters, sample size or position of $k^0_m$ after $m$. For each test statistic, the algorithm steps are given to calculate the corresponding critical values. Afterward, details are given how to calculate empirical test size, empirical test power and to estimate the change-point location. \\
All simulations were performed using the R language. The program codes can be requested from the author.
\subsection{Test by weighted  CUSUM, without bootstrapping}
Firstly, following simulation steps are realized in order to calculate the critical values $c_\alpha(\gamma)$ in accordance with the Corollary \ref{Corollary 1}:
\begin{enumerate}
\item Calculate ${\cal D}\equiv\cro{\A^t \B^{-1}\A }^{1/2}$. 
\item Simulate $M$ replications of the random variable $V_\gamma=\sup_{0 \leq t \leq 1/{\cal D}^2} \frac{(1+t-{\cal D}^2t)|W(t)| }{t^\gamma}$, with $\acc{W(t), 0 \leq t \leq 1/{\cal D}^2}$ a Wiener process, or $V_\gamma=\sup_{0 \leq t \leq  \frac{T}{1+{\cal D}^2T} }\frac{(1+t-{\cal D}^2t)|W(t)|}{t^\gamma}$, with $\acc{W(t), 0 \leq t \leq T/{1+\cal D}^2T}$ a Wiener process, respectively, taking into account the two possible cases \textit{(i)} or \textit{(ii)} concerning $T_m$ of Theorem \ref{theorem 2.1}. 
\item On the basis of  $M$ replications of $V_\gamma$ we  calculate the critical values $c_\alpha(\gamma)$ such that
$\eP[V_\gamma >c_\alpha(\gamma) ]=\alpha$. 
\end{enumerate}
A Brownian motion is generated using  the \textit{BM} function in  R package(\textit{sde}). Once the critical values $c_\alpha(\gamma)$ are available, the change absence against the change of the model is tested using the  statistic $Z_\gamma(m)$, given by relation (\ref{stat_test1}).  In order to calculate the empirical   test size, an without change-point model is considered and we count, the number of times, on the Monte-Carlo replications, when we obtain  $Z_\gamma(m)>c_\alpha(\gamma)$. For the   calculation of  the empirical test power, the hypothesis $H_1$ is considered true,  that there exists a change-point. We fix $T_m=500$, $k^0_m=25$ (or $k^0_m=2$) and  we vary the sample size $m=25, 100, 300$, $\gamma=0,0.25,0.45,0.49$, $\alpha=0.025, 0.05, 0.10$.  For every  combination, $1000$ Monte-Carlo replications are realized. On the 1000 replications, we computed the frequency among which the test statistic  $Z_\gamma(m)$  exceeds the critical value $c_\alpha(\gamma)$. In order to estimate the change-point location, we find the first point $k$ in the interval $1, \cdots, T_m$ such that $(\hat \sigma_m)^{-1} \left| \sum^{m+k}_{i=m+1} \hat \varepsilon_i \right|/g(m,k,\gamma)$  exceeds  critical value $c_\alpha(\gamma)$.\\
For both models, in order to study the importance that $\eff(\x,\eb^0_m)$ is  bounded or not, two regression parameters $\eb^0_m$ after the change-point are considered: one for which $\eff(\x,\eb^0_m)$ is bounded and another for which  $\eff(\x,\eb^0_m)$ is not bounded. Even though the theoretical results are valid, we will study the precision of the change-point location estimator.

\subsubsection{Growth model}
Let us consider first the  growth function  $f(x;\eb)=b_1-\exp(-b_2 x)$ which models many phenomena, with the parameters  $\eb=(b_1,b_2) \in \Theta$, $\Theta \subseteq \R \times \R_+$ compact  and $x \in \R$. In this case the dimension of $\eb$ is 2 ($q=2$) and it there is a single regressor ($p=1$).  We generate the response variable  $X \sim {\cal N}(0,\sigma^2_X)$  and the errors  $\varepsilon \sim {\cal N}(0,0.5)$. The true values of regression parameters before the  change-point are  $\ebo=(0.5, 1)$ and after  $\eb^0_m=(1, 2)$.  By elementary calculations we obtain 
$\eE[X \exp(-b_2 X)]=-b_2 \sigma^2_X \exp( b^2_2 \sigma^2_X/2)$, $\eE[X^2 \exp(-2b_2 X)]=\sigma^2_X[1+4 b^2_2]\exp(2 b^2_2 \sigma^2_X)$, then
\[
\A=\eE[\ef(X,\eb)]=\left[
\begin{array}{c}
1 \\
-b_2 \sigma^2_X \exp( b^2_2 \sigma^2_X/2)
\end{array}
\right]
\]
\[
\B=\eE[\ef(X;\eb)\ef^t(X;\eb)]=\left[
\begin{array}{ccc}
1 & & -b_2 \sigma^2_X \exp( b^2_2 \sigma^2_X/2)\\
& & \\
-b_2 \sigma^2_X \exp(b^2_2 \sigma^2_X/2) & & \sigma^2_X[1+4 b^2_2]\exp(2 b^2_2 \sigma^2_X)
\end{array}
\right].
\]
Obviously ${\cal D}=1$ for any value of  $\sigma^2_X$ and of the parameters $b_1,b_2$. This means that we obtain the same quantiles that in the paper of the Horv\'ath et al.(2004). The empirical quantiles (critical values) $c_\alpha(\gamma)$ of the random variable $V_\gamma$ are given in the Table  \ref{fractiles}.
\begin{table}
 \caption{\footnotesize The $(1-\alpha)$ quantiles (critical values) $c_\alpha(\gamma)$ of the random variable $V_\gamma$ (specified in subsection 5.1, Step 2) calculated on  $50000$ Monte-Carlo replications. Growth model.}
\begin{center}
{\scriptsize
\begin{tabular}{|c|ccccc|}\hline  
$\gamma \downarrow ; \alpha \rightarrow$ & 0.01 & 0.025 & 0.05 & 0.10 & 0.25 \\ \hline 
0 & 2.7959 & 2.5033  & 2.2411 & 1.9595 & 1.5322 \\
0.15& 2.8581 & 2.5690  & 2.3058 & 2.0313 & 1.6146 \\
0.25& 2.9243 & 2.6368  & 2.3841 & 2.1082 & 1.7014\\
0.35& 3.0220 & 2.7536  & 2.5044 & 2.2414 & 1.8462 \\
0.45& 3.2578 & 3.0051  & 2.7878 & 2.5391 & 2.1639 \\
0.49& 3.5214 & 3.2668  & 3.0473 & 2.8040 & 2.4133 \\ \hline 
\end{tabular} 
}
\end{center}
\label{fractiles} 
\end{table}
\begin{table}
 \caption{\footnotesize Empirical sizes of test based on the statistic (\ref{stat_test1}) for a growth model. Calculated for  $1000$ Monte-Carlo replications and $T_m=500$.  }
\begin{center}
{\scriptsize
\begin{tabular}{|c||ccc|ccc|ccc|}\hline
  & $\alpha$= &0.025  &  &$\alpha=$  & 0.05 &  & $\alpha=$ &0.10  &  \\ \hline  
$\gamma \downarrow$& m=25  & m=100 & m=300  & m=25  & m=100 & m=300 & m=25  & m=100 & m=300 \\ \hline 
0 & 0.0051 & 0.0026 & 0.0003 & 0.0075 & 0.0038 & 0.0006 & 0.0132 & 0.0081 & 0.0020 \\
0.25 & 0.0058& 0.0025 & 0.0023 & 0.0084 & 0.0043 & 0.0043 & 0.0154 & 0.0075 & 0.0079 \\
0.45 & 0.0066 & 0.0033 & 0.0023 & 0.0089 & 0.0050 & 0.0043 & 0.0130 & 0.089 & 0.0079 \\
0.49 & 0.0046 & 0.0021 & 0.0014 & 0.0065 & 0.0032 & 0.0026 & 0.0093 & 0.0064 & 0.0070 \\ \hline
\end{tabular} 
}
\end{center}
\label{Risque1} 
\end{table}
Based on these empirical quantiles, we are going to study the test size and its power for various values of  $m$, $\gamma$ and $\alpha$. We realize $1000$  Monte-Carlo replication of the model and we take $T_m=500$. The  empirical test sizes  are presented in Table \ref{Risque1}.  We observe that the obtained values are smaller widely to the fixed  $\alpha$ theoretical size .
 On the 1000 replications we found that empirical test power is 1, in any case.  For the same parameters, we estimate  now  as follows the change-point location.  For $\gamma=0.49, 0.25, \gamma=0$ and $m=25$ or $100$, after $10000$  Monte-Carlo model replications  in Table  \ref{est_cp1} are given the minimum, median, mean, third quartile and maximum of the change-point location estimations. For $m=300$, the  results are similar to those obtained for $m=100$, thus we don't present them. We observe that the obtained change-point estimates are biased, and that  considering either the  median or the mean, there is a delay time in change-point detection.  In the Table \ref{est_cp3} we have the summarized  results when the  change-point is immediately later  after $m$, for $k^0_m=2$. From these two Tables \ref{est_cp1} and \ref{est_cp3} we deduce that, with respect to $\gamma$,  when the change is in  $k^0_m=25$, there is no difference concerning the location change-point precision. If the change is immediately ($k^0_m=2$), the precision decreases when $\gamma$ decreases. \\ 
In  all tables, we indicated between "()" the obtained  results when $\eb^0_m=(1,-0.5)$, case in which the function $\eff(\textbf{x};\eb^0_m)$ is not bounded for all  $\textbf{x}$. The results are worse, even though the break in  $k^0_m$ is largest.   
\begin{table}
 \caption{\footnotesize Estimation of the change-point location based on the statistic (\ref{stat_test1}), for  $10000$ Monte-Carlo replications, $T_m=500$,  $k^0_m=25$, $\ebo=(0.5, 1)$, $\eb^0_m=(1, 2)$ and between  () for $\eb^0_m=(1,-0.5)$.   Growth model.}
\begin{center}
{\scriptsize
\begin{tabular}{|c||c||ccc|ccc|}\hline
&   & $m$= &25  &  &$m=$  & 100 &   \\ \hline  
$\gamma$ & $summary(\hat \tau_m) \downarrow$ ; $\alpha \rightarrow$& $0.025$  & $0.05$ & $0.10$  & $0.025$  & $0.05$ & $0.10$ \\ \hline 
0.49 & \textit{min} & 1 (1) &  1 (1) & 1 (1) &1  (1) &1  (1) & 1 (1)  \\
& \textit{median(Q2)} &32 (41)&32 (39) & 31 (37) &31  (37) &31  (35) &31  (34)  \\
& \textit{mean} &35 (46) & 34 (43) &33 (40) &34  (39) &33  (38) &32 (36)  \\
& \textit{Q3} &39 (54)&38 (51) &37 (48) &37 (45) &37 (43) &36 (41)  \\
& \textit{max} &148 (279) &148 (248) &148 (232) &112 (144) &111 (143) &109 (132) \\ \hline \hline
0.25 & \textit{min} &1 (1) &1 (1) &1 (1) &9 (1) & 9 (1) &9 (1) \\
& \textit{median(Q2)} &32 (39)&31 (37) &31 (35) & 32 (38) & 32 (37) &31 (35)  \\
& \textit{mean} &34 (43) &33 (40) &32 (37) & 34 (41) &34  (39) &33 (38)  \\
& \textit{Q3} &39 (49)&38 (46) &36 (43) &38 (46) &37 (44) &37 (42) \\
 & \textit{max} &141 (245) &141 (219) &119 (200) &110 (137) & 110 (125) &95 (125) \\ \hline \hline
 0 & \textit{min} & 1 (1) & 1 (1) & 1 (1) & 9 (25) &9  (21) & 9 (20) \\
& \textit{median(Q2)} &33  (41)& 32 (39) & 31  (36) & 34  (43) & 33 (41) & 33 (39)  \\
 & \textit{mean} & 35  (45) & 34  (42) & 33  (39) & 36  (45) & 35  (43) & 35 (41)  \\
& \textit{Q3} & 39  (52)& 38 (49) & 37 (45) &41  (52) &  40 (49) & 39 (46)  \\
& \textit{max} & 121  (300) & 115 (284) & 115  (223) & 115 (136) & 115 (132) & 115 (131)  \\ \hline \hline 
\end{tabular} 
}
\end{center}
\label{est_cp1} 
\end{table}

\begin{table}
 \caption{\footnotesize Estimation of the change-point location based on the statistic (\ref{stat_test1}), for  $10000$ Monte-Carlo replications, $T_m=500$,  $k^0_m=2$, $\ebo=(0.5, 1)$, $\eb^0_m=(1, 2)$ and between () for $\eb^0_m=(1,-0.5)$.  Growth model.}
\begin{center}
{\scriptsize
\begin{tabular}{|c||c||ccc|ccc|}\hline
 &  & $m$= &25  &  &$m=$  & 100 & \\ \hline  
$\gamma $ &$summary(\hat \tau_m)  \downarrow$; $\alpha \rightarrow$& $0.025$  & $0.05$ & $0.10$  & $0.025$  & $0.05$ & $0.10$ \\ \hline 
0.49 & \textit{min} & 1 (1) &  1 (1) & 1 (1) &1  (1) & 1 (1) & 1 (1)  \\
& \textit{median(Q2)} & 6 (6)& 6 (6) & 5 (5) & 6 (6) & 6 (6) & 5 (5)  \\
& \textit{mean} & 8 (10) & 8 (9) & 7 (8) & 8 (9) & 7 (8) & 7 (7)  \\
& \textit{Q3} & 10 (12)& 10 (11) & 9 (10) & 10 (11) & 9 (10) & 9 (9)  \\
& \textit{max} & 109 (224) & 109 (218) & 79 (185) & 91 (81) & 91 (71) & 91 (67)  \\ \hline \hline
0.25 & \textit{min} & 1 (1) &1  (1) & 1 (1) & 1 (2) & 1 (1) & 1 (1) \\
& \textit{median(Q2)} & 7 (8)& 6 (7) & 6 (6) & 7 (10) & 7 (9) & 7 (8)  \\
& \textit{mean} & 9 (11) & 8 (10) & 7 (9) & 10 (12) & 9 (11) & 9 (10) \\
& \textit{Q3} & 11 (14)& 10 (12) & 10 (11) & 13 (16) & 12 (15)& 11 (13)  \\
& \textit{max} & 74 (156) & 72 (133) & 72 (133) & 89 (94) & 87 (94) & 85 (94)  \\ \hline \hline 
0 & \textit{min} &1 (1) & 1 (1) & 1 (1) & 6 (3) & 6 (3) & 5 (3)  \\
& \textit{median(Q2)} & 8 (11)& 7 (10) & 7 (9) & 10 (17) & 9 (15) & 9 (14)  \\
& \textit{mean} & 10 (14) & 9 (13) & 9 (11) & 12 (19) & 12 (17) & 11 (15)  \\
& \textit{Q3} & 13 (18)& 12 (16) & 11 (14) & 16 (25) & 15 (23) & 15 (20)  \\
& \textit{max} & 93 (210) & 88 (179) & 88 (175) & 91 (131) & 76 (124) & 73 (111)  \\ \hline
\end{tabular} 
}
\end{center}
\label{est_cp3} 
\end{table}
\subsubsection{Compartmental model}
Another very interesting nonlinear model, with numerous applications, is the compartmental model. Examples and references of important applications for these models are given in Seber and Wild(2003) (see also the references therein): it describes  the movement of lead in the human body, the kinetics of drug movement when the drug is injected at an intramuscular site, etc... Consider two-compartment function $h_\beta(x)=b_1 \exp(-b_1 x)+b_2 \exp(-b_2x)$, $\eb=(b_1,b_2) \in \Theta \subseteq \R^2_+$. In this case $q=2$ and $p=1$.\\
As for the growth example,  we consider a gaussian  response variable $X \sim {\cal N}(0, \sigma^2_X)$. For this model we have $\eE[\exp(-bX)]=\exp(b^2 \sigma^2_X/2)$, $\eE[X \exp(-bX)]=-b \sigma^2_X \exp(b^2 \sigma^2_X/2)$, $\eE[X^2 \exp(-2b X)]=\sigma^2_X[1+4 b^2]\exp(2 b^2 \sigma^2_X)$. Then
\[
\A=\eE[\ef(X,\eb)]=\left[
\begin{array}{c}
(1+b_1 \sigma^2_X)\exp(b^2_1 \sigma^2_X/2) \\
-(1+b_2 \sigma^2_X)\exp(b^2_2 \sigma^2_X/2)
\end{array}
\right].
\]
And with the notations $B_{11}=1+b_1^2 \sigma^2_X(5+4b_1^2)\exp (2 b_1^2 \sigma^2_X)$, $B_{12}=1+\sigma^2_X [ (b_1+b_2)^2+b_1b_2(1+(b_1+b_2)^2)] \exp ((b_1+b_2)^2 \sigma^2_X/2)$, $B_{22}=1+b_2^2 \sigma^2_X(5+4b_2^2)\exp (2 b_2^2 \sigma^2_X)$, we have the matrix 
 \[
\B=\eE[\ef(X;\eb)\ef^t(X;\eb)]
=\left[
\begin{array}{ccc}
 B_{11}& & B_{12}\\
& & \\
B_{12} & & B_{22}
\end{array}
\right].
\]
Contrary to the previous case, the value of ${\cal D}$ depends on the variance $\sigma^2_X$ of the random variable $X$ and on the  parameters of the growth function. Hence, for each value of  $\ebo$ and of  variance of $X$ we need to calculate the quantiles. For the simulations, let us consider $\sigma^2_X=1$  and $\ebo=(1.2,1)$. In this case ${\cal D}=0.5741$.\\
The empirical quantiles $c_\alpha(\gamma)$ of the random variable $V_\gamma$, specified at the beginning of this subsection, are given in the Table \ref{fractiles_comp}.\\
The simulations are carried out for historical data of size $m=25, 100 $ or $300$ and $T_m=500$ observation after $m$.  The empirical type I error probabilities are presented in the Table \ref{Risque1_comp} calculated by 1000 Monte-Carlo replications.
As for the growth example, the empirical  power test is 1 for each  value of  $\gamma,\alpha$, when $k^0_m=25$ and $\eb^0_m=(1,2)$.\\
 In Tables \ref{est_cp1_comp} and \ref{est_cp3_comp}, the summarized results on the change-point estimations obtained on 10000 Monte-Carlo replications, varying $m$, $\gamma$ and theoretic test size $\alpha$.   Between "()" we give the results for  $\eb^0_m=(-0.5,2)$,  when the function $\eff(\textbf{x};\eb^0_m)$ is not bounded for all value of $\textbf{x}$.\\
We can make the following observations. As for the growth example, the results are less good in the  case $\eff(\textbf{x};\eb)$ not bounded: the method  detects later the change and especially we have   greater maximal values for the change-point estimation  $\hat \tau_m$. In the two case, $k^0_m=25$ and $k^0_m=2$, the  precision of  $\hat \tau_m$ decreases when $\gamma$ decreases. The change-point estimation is more precise than for the growth model.
\begin{table}
 \caption{\footnotesize The $(1-\alpha)$ quantiles (critical values) $c_\alpha(\gamma)$ of the random variable $V_\gamma$ (specified in subsection 5.1, Step 2) calculated on  $50000$ Monte-Carlo replications. Compartmental model,  $\ebo=(1.2,,1)$, $\sigma^2_X=1$. }
 \begin{center}
{\scriptsize
\begin{tabular}{|c|ccccc|}\hline  
$\gamma \downarrow$ ; $\alpha \rightarrow$ & 0.01 & 0.025 & 0.05 & 0.10 & 0.25 \\ \hline 
0 & 6.2165 & 5.5233  & 4.9211 & 4.2812 & 3.2689 \\
0.15& 5.7627 & 5.1279  & 4.5862 & 4.0014 & 3.0854 \\
0.25& 5.4929 & 4.9022  & 4.3838 & 3.8395 & 2.9833\\
0.35& 5.2355 & 4.6960  & 4.2092 & 3.7024 & 2.9142 \\
0.45& 5.0383 & 4.5223  & 4.0786 & 3.5998 & 2.9191 \\
0.49& 4.9682 & 4.4702  & 4.0555 & 3.6032 & 2.9945 \\ \hline 
\end{tabular} 
}
\end{center}
\label{fractiles_comp} 
\end{table}

\begin{table}
 \caption{\footnotesize Empirical sizes of test based on the statistic (\ref{stat_test1}) for a compartmental model. Calculated for  $1000$ Monte-Carlo replications and $T_m=500$. }
 \begin{center}
{\scriptsize
\begin{tabular}{|c||ccc|ccc|ccc|}\hline
  & $\alpha$= &0.025  &  &$\alpha=$  & 0.05 &  & $\alpha=$ &0.10  &  \\ \hline  
$\gamma$& m=25  & m=100 & m=300  & m=25  & m=100 & m=300 & m=25  & m=100 & m=300 \\ \hline 
0 & 0.0003 & 0 &0 & 0.0005 & 0 & 0 & 0.0007 & 0 & 0 \\
0.25 & 0.0006 & 0 & 0 & 0.0006 &  0& 0 & 0.0009 & 0 & 0 \\
0.45 & 0.0010 & 0 & 0 & 0.0013 & 0 & 0 & 0.0014 & 0.0002 & 0 \\
0.49 & 0.0009 & 0 & 0 & 0.0011 & 0.0001 & 0.0001 & 0.0012 & 0.0004 &  0.0002 \\ \hline
\end{tabular} 
}
\end{center}
\label{Risque1_comp} 
\end{table}

\begin{table}
 \caption{\footnotesize Estimation of the change-point location based on the statistic (\ref{stat_test1}), for  $10000$ Monte-Carlo replications, $T_m=500$, $k^0_m=25$, $\ebo=(1.2, 1)$, $\eb^0_m=(1, 2)$ and between () for $\eb^0_m=(-0.5,2)$. Compartmental model.}
\begin{center}
{\scriptsize
\begin{tabular}{|c||c||ccc|ccc|}\hline
&   & $m$= &25  &  &$m=$  & 100 &   \\ \hline  
$\gamma $ & $summary(\hat \tau_m)\downarrow$ ; $\alpha \rightarrow$& $0.025$  & $0.05$ & $0.10$  & $0.025$  & $0.05$ & $0.10$ \\ \hline 
0.49 & \textit{min} & 1 (1) & 1 (1) & 1 (1) & 26 (1) & 5 (1) & 1 (1)  \\
& \textit{median(Q2)} & 30 (31)& 30 (31) & 29 (31) & 29 (31) & 29 (30) & 29 (30)  \\
& \textit{mean} & 30 (34) & 30 (33) & 30 (33) & 30 (32) & 30 (32) & 29 (32)  \\
& \textit{Q3} & 33 (37)& 33 (37) & 32 (36) & 32 (36) & 32 (35) & 31 (34)  \\
 & \textit{max} & 63 (105) & 63 (102) & 60 (102) & 55 (105) & 53 (105) & 53 (102) \\ \hline \hline 
0.25 & \textit{min} & 1 (1) & 1 (1) & 1 (1) & 26 (26) & 26 (26) & 26 (26)  \\
& \textit{median(Q2)} & 31 (32)& 30 (32) & 30 (31) & 30 (32) & 30 (32) & 30 (31)  \\
& \textit{mean} & 31 (35) & 31 (34) & 30 (34) & 31 (34) & 31 (34) & 30 (33) \\
& \textit{Q3} & 34 (39)& 34 (38) & 33 (37) & 34 (38) & 34 (38) & 33 (37)  \\
& \textit{max} & 69 (124) & 69 (124) & 59 (115) & 61 (99) & 60 (99) & 52 (97)  \\ \hline \hline 
0 & \textit{min} & 3 (1) & 3 (1) & 3 (1) & 26 (26) & 26 (26) &  26 (26)  \\
& \textit{median(Q2)} & 32 (34)& 31 (33) & 31 (33) & 32 (35) & 32 (34) & 31 (33)  \\
& \textit{mean} & 33 (37) & 32 (36) & 32 (35) & 33 (37) & 33 (36) & 32 (35)  \\
& \textit{Q3} & 36 (41)& 35 (40) & 34 (39) & 37 (42) & 36  (41) & 35 (40)  \\
& \textit{max} & 73 (143) & 67 (138) & 63 (114) & 69 (120) & 69 (120) & 66 (102)  \\ \hline \hline 
\end{tabular} 
}
\end{center}
\label{est_cp1_comp} 
\end{table}

\begin{table}
 \caption{\footnotesize Estimation of the change-point location based on the statistic (\ref{stat_test1}), for  $10000$ Monte-Carlo replications, $T_m=500$,  $k^0_m=2$, $\ebo=(1.2, 1)$, $\eb^0_m=(1, 2)$ and between () for $\eb^0_m=(-0.5,2)$. Compartmental model.}
\begin{center}
{\scriptsize
\begin{tabular}{|c||c||ccc|ccc|}\hline
&   & $m$= &25  &  &$m=$  & 100 &    \\ \hline  
$\gamma$ & $summary(\hat \tau_m)\downarrow$ ; $\alpha \rightarrow$& $0.025$  & $0.05$ & $0.10$  & $0.025$  & $0.05$ & $0.10$\\ \hline 
 & \textit{min} &  1 (1) & 1 (1) & 1 (1) & 3 (3) & 3 (2) & 3 (1)  \\
0.49 & \textit{median(Q2)} & 5 (6)& 4 (5) & 4 (5) &  5 (5)& 4 (5) & 4 (5) \\
& \textit{mean} & 5 (7)& 5 (7) & 5 (6) & 5 (7) & 5 (6) & 5 (6) \\
& \textit{Q3} & 7 (9)& 6 (8) & 6 (8) & 7 (8) & 6 (8) & 6 (7)  \\
& \textit{max} & 29 (68) & 29 (61) & 29 (60) & 27 (80) & 27 (80) & 27 (49)  \\ \hline \hline 
0.25 & \textit{min} & 1 (1) & 1 (1) & 1 (1) & 3 (3) &  3 (3)& 3 (3)  \\
& \textit{median(Q2)} & 5 (7)& 5 (6) & 5 (6) & 6 (7)& 6 (7) & 5 (7)  \\
& \textit{mean} & 6 (9) & 6 (8) & 6 (8) & 7 (9) & 6 (9) & 6 (9) \\
& \textit{Q3} & 8 (11)& 8 (11) & 7 (10)& 9 (12) & 9 (12) & 8 (11)  \\
& \textit{max} & 37 (86) & 37 (86) & 35 (86) & 36 (98) & 35 (97) & 34 (90)  \\ \hline \hline 
0 & \textit{min} & 1 (1) & 1 (1) & 1 (1) & 3 (3) & 3 (3) & 3 (3)  \\
& \textit{median(Q2)} & 7 (8)& 6 (8) & 6 (8) & 9 (11)& 9 (10) & 8 (9)  \\
& \textit{mean} & 7 (11) & 7 (10) & 7 (10) & 9 (13) & 9 (12) & 8 (12)  \\
& \textit{Q3} & 10 (14)& 9 (14) & 9 (13)& 13 (16) & 12 (15) & 11 (15)  \\
& \textit{max} & 35 (99) & 34 (91) & 34 (91) & 49 (110) & 40 (91) & 40 (91)  \\ \hline \hline 
\end{tabular} 
}
\end{center}
\label{est_cp3_comp} 
\end{table}

\subsection{Test using the bootstrapping}
In this case, the  calculation of the critical values  $c_{m,k;\alpha}(\gamma)$ defined by (\ref{eq19}) is more laborious. We go to see if the simulation results are better than by weighted CUSUM without bootstrapping, case in which it deserves to make calculation effort.   \\
We now describe in detail the algorithm steps for  calculate the critical values $c_{m,k;\alpha}(\gamma)$.\\
\textbf{Step 1.} We fix $\alpha, \gamma, N, L, m, T_m$ (see the notations given in Section 4 for N and L).\\
\textbf{Step 2.} \begin{itemize}
\item We calculate $J=T_m/L$;
\item For $j=0,1, \cdots, (J-1)L$, the following  random variable  are generated  
\[
\tilde V_j \equiv \frac{1}{\hat \sigma^{(*)}_{m,j}} \sup_{1 \leq l \leq T_m} | \tilde \Gamma(m,j,l,\gamma)( \varepsilon^*_{m,j}(1), \cdots, \varepsilon^*_{m,j}(m+l) )|
\]
\end{itemize}
\textbf{Step 3.} For $\tilde j=0, 1, \cdots, (J-1)L$, we generate the random variables  $\tilde W_{\tilde j}$ which are mixtures of the random variables $\tilde V_j$ generated to step  2. \\
For each $\tilde j=1, \cdots, (J-1)L$, we generate  a multinomial distribution with parameters 1(number of trials) and the probability vector $p_{\tilde j}=( 1/\tilde j, \cdots, 1/\tilde j)$. On the basis of this, thus, $\tilde W_{\tilde j}=\tilde V_j$ for $j=0,1, \cdots, \tilde j-1$ with the  probability $1/\tilde j$. \\
 \textbf{Step 4.} We repeat the steps  2 and 3 making  $M$ Monte-Carlo replications. At the end,  we shall have $M$ realizations for every random variable $\tilde W_{\tilde j}$, $\tilde j=0, 1, \cdots, (J-1)L$. \\
 \textbf{Step 5.} We calculate for every $k=jL, jL+1, \cdots, (j+1)L-1$ for $j=0,1, \cdots, J$ the random variables $\tilde U_k=\tilde W_j$.\\
 \textbf{Step 6.} On the basis of $M$ replications, for each $k=1, \cdots, T_m$, we calculate the critical values  $c_{m,k;\alpha}(\gamma)$ such that $\eP[\tilde U_k >c_{m,k;\alpha}(\gamma)]=\alpha$.\\ 
 
 The change absence against the change of the model is tested using the statistic $Z^{(b)}_{\gamma;\alpha}(m)$ given by (\ref{stat_test2}). In  order to calculate the empirical   test size, an without change-point model is considered and we count, the number of times, on the  Monte Carlo replications, when we obtain $Z^{(b)}_{\gamma;\alpha}(m)>1$. Recall that the change-point estimation  $\hat \tau^{(b)}_m$ is calculated using relation (\ref{eq7_b}). \\
 Let us consider  $m=25$ and $m=100$. For $m=300$, the  results are similar to those obtained for $m=100$, thus we don't present them. In the  case $m=100$ we consider $L=m/50$ and in the case $m=25$ we take $L=m/10$. For $\gamma$ we take only two values: 0.25 et 0.49. If $k^0_m=25$, the empirical power test is 1 in all cases:  for the two model type (growth or compartmental) and for the every parameters $\gamma$ and $k^0_m$.\\
 The same parameter settings are used as in the previous simulation study, in the subsection 5.1.

\subsubsection{Compartmental model}
 The empirical test size based on the statistic $\tilde U_k$ (of Step 5), calculated for 1000 Monte-Carlo replications and $T_m=500$,  are given in the Table  \ref{Risque1_comp_boot}. By comparing the  Tables \ref{Risque1_comp} and  \ref{Risque1_comp_boot}, we deduce that the empirical test sizes  are smaller by the  bootstrapping method. \\
 The results concerning  $\hat \tau^{(b)}_m$, the estimation of $k^0_m$, presented in the Tables \ref{est_cp1_comp_boot} and \ref{est_cp3_comp_boot}, are almost the same for  $\gamma=0.49$ and $\gamma=0.25$. Apart from $\gamma=0.25$ and $k^0_m=25$, the  results for $\hat \tau^{(b)}_m$ are not better than those obtained by the method without bootstrapping. 
 
\begin{table}
 \caption{\footnotesize Empirical sizes of test based on the statistic $\tilde U_k$ given in subsection 5.2, Step 5, for a compartmental model, for bootstrap critical values. Calculated for  $1000$ Monte-Carlo replications and $T_m=500$.  }
\begin{center}
{\scriptsize
\begin{tabular}{|c||cc|cc|cc|}\hline
  & $\alpha$= &0.025    &$\alpha=$  & 0.05   & $\alpha=$ 0.10  &  \\ \hline  
$\gamma$& m=25  & m=100   & m=25  & m=100  & m=25  & m=100  \\ \hline 
0.25 & 0.0001 & 0  & 0.0003 &  0 & 0.0007 & 0  \\
0.49 & 0.0004 & 0 & 0.0006 & 0 & 0.0008 & 0  \\ \hline
\end{tabular} 
}
\end{center}
\label{Risque1_comp_boot} 
\end{table}
\begin{table}
 \caption{\footnotesize Estimation of the change-point location based on the statistic (\ref{stat_test2}), for  $10000$ Monte-Carlo replications, $T_m=500$,  $k^0_m=25$, $\ebo=(1.2, 1)$, $\eb^0_m=(1, 2)$ and between () for $\eb^0_m=(-0.5,2)$. Compartmental model.}
\begin{center}
{\scriptsize
\begin{tabular}{|c||c||ccc|ccc|}\hline
 &  & $m$= &25  &  &$m=$  & 100 &   \\ \hline  
$\gamma$ & $summary(\hat \tau^{(b)}_m)\downarrow$; $\alpha \rightarrow$& $0.025$  & $0.05$ & $0.10$  & $0.025$  & $0.05$ & $0.10$ \\ \hline 
0.49 & \textit{min} & 2 (2) & 1 (1) & 1 (1) & 26 (26) & 26 (26) & 14 (26)  \\
& \textit{median(Q2)} & 32 (30)& 30 (28) & 26 (27) & 32 (35) & 31 (33) & 29 (32)  \\
& \textit{mean} & 32 (33) & 30 (38) &  27 (29) &  32  (37) & 31 (36) & 30 (34) \\
& \textit{Q3} & 36 (36)& 33 (35) & 30 (33) & 35 (42) & 34 (40) & 32 (38)  \\
& \textit{max} & 90 (118)& 78 (118) & 65 (92) & 72 (122) & 63 (122) & 54 (122)  \\ \hline \hline 
0.25 & \textit{min} & 5 (1) & 5 (1) & 1 (1) &  26 (26) & 26 (26) & 26 (26)  \\
& \textit{median(Q2)} & 28 (30)& 26 (28) & 24 (26) & 32 (35) & 31 (34) & 30 (32)  \\
& \textit{mean} & 28 (32) & 27 (30)& 25 (28) & 32 (38) & 32 (36) & 31 (34) \\
& \textit{Q3} & 32 (38)& 28 (34) & 28 (31) & 34 (41) & 34 40) & 34 (38)  \\
& \textit{max} & 66 (106) & 59 (103) & 57 (103) & 71 (134) & 66 (118) & 57 (98)  \\ \hline
\end{tabular} 
}
\end{center}
\label{est_cp1_comp_boot} 
\end{table}

\begin{table}
 \caption{\footnotesize Estimation of the change-point location based on the statistic (\ref{stat_test2}), for  $10000$ Monte-Carlo replications,  $T_m=500$,  $k^0_m=2$, $\ebo=(1.2, 1)$, $\eb^0_m=(1, 2)$and between () for $\eb^0_m=(-0.5,2)$. Compartmental model.}
\begin{center}
{\scriptsize
\begin{tabular}{|c||c||ccc|ccc|}\hline
&   & $m$= &25  &  &$m=$  & 100 &   \\ \hline  
$\gamma=0.49$ & $summary(\hat \tau^{(b)}_m)\downarrow$; $\alpha \rightarrow$& $0.025$  & $0.05$ & $0.10$  & $0.025$  & $0.05$ & $0.10$ \\ \hline 
0.49 & \textit{min} & 3 (1) & 2 (1) & 2 (1) & 3 (3) & 3 (3) & 3 (3) \\
 & \textit{median(Q2)} & 6 (8)& 6 (6) & 5 (5) & 6  (7)& 5 (7) & 5 (6)  \\
& \textit{mean} & 8 (11)& 7 (9) &  6 (7) & 7 (10) & 6 (10) & 5 (8)  \\
& \textit{Q3} & 10 (13)& 9 (11) & 8 (9) & 9 (14) & 7  (12) & 7 (10)  \\
& \textit{max} & 65 (76) & 64  (76) & 33 (63) & 43  (86) &  42 (75) & 27 (75)  \\ \hline \hline 
0.25 & \textit{min} & 3 (2) & 1 (1) & 1 (1) & 3 (3) & 3 (3)& 3 (3 ) \\
& \textit{median(Q2)} & 6 (7)& 5 (6) & 4 (5) & 7 (10)& 6 (10) & 6 (7) \\
& \textit{mean} & 7 (10) & 6 (9) & 5 (7) & 8 (13) & 7 (12) & 7 (9) \\
& \textit{Q3} & 9 (13)& 8 (12) & 7 (9) & 11 (18) & 10 (16) &  9 (12)  \\
& \textit{max} & 38  (81) & 32 (81)& 28 (62) & 37 (119) & 36 (74) & 32 (74) \\ \hline \hline
\end{tabular} 
}
\end{center}
\label{est_cp3_comp_boot} 
\end{table}

\subsubsection{Growth model}
Tables \ref{Risque1} and \ref{Risque1_growth_boot} indicate that the empirical test size obtained using the bootstrapped critical values are sharply lower than empirical test size without bootstrapping. Concerning the change-point estimation (Table \ref{est_cp1_growth_boot}), for $\gamma=0.49$, $m=25$  and $k^0_m=25$, the results for   $\hat \tau^{(b)}_m$ are better than by the weighted  CUSUM method without bootstrapping. On the other hand, for $\gamma=0.49$, $m=100$, the  results are less good using the bootstrapped critical values.

\begin{table}
 \caption{\footnotesize Empirical sizes of test based on the statistic $\tilde U_k$ given in subsection 5.2, Step 5, for a growth model, for bootstrap critical values. Calculated for  $1000$ Monte-Carlo replications and $T_m=500$.  }
\begin{center}
{\scriptsize
\begin{tabular}{|c||cc|cc|cc|}\hline
  & $\alpha$= &0.025    &$\alpha=$  & 0.05   & $\alpha=$ 0.10  &  \\ \hline  
$\gamma$& m=25  & m=100   & m=25  & m=100  & m=25  & m=100  \\ \hline 
0.25 & 0.0005 & 0 & 0.0012 & 0 & 0.0027 & 0.0004  \\
0.49 & 0.0002 & 0 & 0.0003 & 0 & 0.0007 & 0  \\ \hline
\end{tabular} 
}
\end{center}
\label{Risque1_growth_boot} 
\end{table}

\begin{table}
 \caption{\footnotesize Estimation of the change-point location based on the statistic (\ref{stat_test2}), for  $10000$ Monte-Carlo replications,  $T_m=500$, $\gamma=0.49$,  $\ebo=(0.5, 1)$, $\eb^0_m=(1, 2)$ and between () for $\eb^0_m=(1,-0.5)$. Growth model.}
\begin{center}
{\scriptsize
\begin{tabular}{|c||c||ccc|ccc|}\hline
&  & $m$= &25  &  &$m=$  & 100 &  \\ \hline  
$k^0_m$ & $summary(\hat \tau^{(b)}_m)\downarrow$ ; $\alpha \rightarrow$& $0.025$  & $0.05$ & $0.10$  & $0.025$  & $0.05$ & $0.10$ \\ \hline 
25 & \textit{min} &  1 (1)&  1 (1) &{ 1 (1)} & {  1 (1)} & {  1 (1)} & {  1 (1)} \\
& \textit{median(Q2)} & {  27 (32)}&  { 26  (30)} &   {25 (28)}  &  {38 (53)} &   {36 (45)} &  {34 (42)}  \\
& \textit{mean} &  {28 (35)} &  {27 (33)} & { 26 (30)} &  {40 (57)} &  {39 (51)} & {  36 (44)}  \\
& \textit{Q3} &  {32 (39)}& { 32  (39) }&  {30 (36) }&  {48  (61)} & {  44 (60) }& {42 (50) } \\
& \textit{max} & { 107 (242)} & { 107  (182)} &  {107  (175)} & { 146  (220)} &  { 145 (184)} &  {140 (156)} \\ \hline \hline 
2 & \textit{min} &  {1 (1) }&  { 1 (1)} & { (1) 1} &  {1 (2)} &  {1 (1)} &  {1 (1) }\\
& \textit{median(Q2)} & { 5 (6)}& {  5 (5)} &  { 5 (5)}  &  {8 (11)} &  {8 (10)} &  {7 (9) } \\
& \textit{mean} &  { 7 (8)} &  {7 (8)} &  {6 (7)} &  {11 (19) }&  {10 (15)} & { 9 (12) } \\
& \textit{Q3} &  { 9 (10)} &  {8 (9)} & { 8 (8)} & { 14 (26)} & { 13 (20)} &  {12 (16)}  \\
& \textit{max} &  {82 (122) } &  {82 (122)} & { 80 (107)} & { 122 (192)} & { 104 (150) }&  {101 (100) }\\ \hline
\end{tabular} 
}
\end{center}
\label{est_cp1_growth_boot} 
\end{table}

\subsection{Conclusion on the simulations}
Two test statistics and their critical regions are, using weighted CUSUM method without and with bootstrapping for two nonlinear models.  In both cases, the empirical sizes are widely smaller than the fixed  theoretical size $\alpha$. But the empirical  sizes of test are without thinking smaller when the critical values are calculated by  bootstrapping. The power test is  equal to 1 for any value of $m$, $\gamma$, $k^0_m$, or theoretic test size $\alpha$. The both test statistics  (\ref{stat_test1}) and (\ref{stat_test2}) detect the change produced in the model. \\
The parameter $\gamma$ does not modify the type I error probability.
Concerning the change-point estimation precision, it does not improve in a significant way by the bootstrapping method or when the number $m$ of historical data increases. This precision can be influenced by $\gamma$ value when the test statistic (\ref{stat_test1}), without bootstrapping, is used.    
It is worth mentioning that the  obtained estimations of  $k^0_m$ by the both methods are slightly biased, the delay time is of order $\simeq +6$ observations, either  for  $m=25$ or for $m=100$ observations. \\
Finally, if  $\eff(\x,\eb)$ is not  bounded, the both test statistics  detect  the change-points, but the estimator bias of  $k^0_m$  increases, if the change is 2 observations after $m$ or 25 observations after $m$. 

\section{Proofs of the Theorems and Propositions}
Here we present the proofs of the results stated in Sections 3 and 4.\\

\noindent {\bf Proof of Theorem  \ref{theorem 2.1}}\\
The proof follows the structure of the Theorem 2.1 proved by Horv\'ath et al. (2004) for the linear case.\\
\textit{(i)} Using  Lemma \ref{Lemma 5.2} and Lemma \ref{Lemma 5.3}   we have
\[
\sup_{1 \leq k <\infty}\sum^{m+k}_{i=m+1} \hat \varepsilon_i/g(m,k,\gamma)  = \sup_{1 \leq k <\infty} \pth{\sum^{m+k}_{i=m+1} \varepsilon_i -\frac{k}{m}\sum^m_{i=1} A_i \varepsilon_i}/g(m,k,\gamma)(1+o_{\eP}(1)) \]
\begin{equation}
\label{F1}
= \sigma \sup_{1 \leq k <\infty} [W_{1,m}(k)-{\cal D} \frac{k}{m} W_2(m)]/g(m,k,\gamma)(1+o_{\eP}(1)).
\end{equation}
with $W_{1,m}$ and $W_2$ two independent Wiener processes   on $[0,\infty)$.
We obtain in a similar way as in the linear case (Theorem 2.1 of Horv\'ath et al., 2004)
\[
\sup_{1 \leq k <\infty} \frac{|W_{1,m}(k) - \frac{k}{m}{\cal D} W_{2,m}(m) |}{g(m,k,\gamma)} \overset{\cal L}{=} \sup_{1 \leq k <\infty} \frac{| W_1(k)- \frac{k}{m}{\cal D} W_2(m)|}{g(m,k,\gamma)},
\]
where $\{W_1(t) \}$, $\{W_2(t) \}$ are two independent Wiener processes on $[0,\infty)$. For all $K>0$, by the 
continuity of $\{W_1(t) -{\cal D}t W_2(1)/(t/(1+t))^\gamma \}$ on $[ 0,K]$ we have
\begin{equation}
\label{LM}
\max_{1 \leq k \leq mK} \frac{| W_1(k)- \frac{k}{m}{\cal D} W_2(m)|}{g(m,k,\gamma)} \overset{\cal L}=\max_{1 \leq k \leq mK} \frac{|W_1 \pth{\frac{k}{m}} -\frac{k}{m}{\cal D} W_2(1) |}{\pth{1+\frac{k}{m}} \pth{\frac{k}{m+k}}^\gamma}
\overset{{a.s.}} {\underset{m \rightarrow \infty}{\longrightarrow}} \sup_{0 \leq t \leq K} \frac{|W_1(t) -{\cal D}t W_2(1) |}{(1+t) \pth{\frac{t}{1+t}}^\gamma}.
\end{equation}
The relations (5.9) and (5.10) of Horv\'ath et al.(2004) hold, then, for all $\delta >0$,
\[ \lim_{K \rightarrow \infty} \limsup_{m \rightarrow \infty} \eP \cro{ \left| \sup_{mK \leq k < \infty} |W_1(\frac{k}{m}) - \frac{k}{m}{\cal D} W_2(1) | /(1+\frac{k}{m})(\frac{k}{m+k})^\gamma -{\cal D} W_2(1) \right|> \delta}=0, 
\]
\[
\lim_{K \rightarrow \infty}\eP \cro{\left| \sup_{K <t<\infty} \frac{|W_1(t)-{\cal D}t W_2(1)|}{(1+t)(\frac{t}{1+t})^\gamma}- {\cal D} W_2(1) \right|> \delta}=0,
\]
thus
\begin{equation}
\label{eq5.14}
\sup_{1 \leq k < \infty }\frac{|W_{1,m}(k) - \frac{k}{m}{\cal D} W_{2,m}(m) |}{g(m,k,\gamma)} \overset{{\cal L}} {\underset{m \rightarrow \infty}{\longrightarrow}} \sup_{0 \leq t < \infty} \frac{|W_1(t)-{\cal D}t W_2(1)|}{(1+t)(\frac{t}{1+t})^\gamma}.
\end{equation}
Let us consider the random  processes $Z(t)=W_1(t)-{\cal D}t W_2(1)$ and $U(t)=(1+{\cal D}^2t) W \pth{\frac{t}{1+{\cal D}^2t}}$, with $\{W(t), 0 \leq t < \infty \}$ a Wiener process. Their variances are 
$Var[Z(t)]=t+{\cal D}^2t^2=t(1+{\cal D}^2t)$, $Var[U(t)]=(1+{\cal D}^2t)^2 t/(1+{\cal D}^2t)=t(1+{\cal D}^2t)$. For $t_1<t_2$, $Cov(Z(t_1),Z(t_2))=\eE[Z(t_1)Z(t_2)]+{\cal D}^2t_1t_2=t_1+{\cal D}^2t_1t_2=t_1(1+{\cal D}^2t_2)$ and since $t/(1+{\cal D}^2t)$ is increasing in $t$, $Cov(U(t_1),U(t_2))=(1+{\cal D}^2t_1)(1+{\cal D}^2t_2) t_1/ (1+{\cal D}^2t_1)=t_1(1+{\cal D}^2t_2)$. Thus, their variances and covariances coincide,  we have $Z(t) \overset{\cal L}{=} U(t)$, for $0 \leq t < \infty$. Let us make the change of variable  $t/(1+{\cal D}^2t)=y$, hence 
\begin{equation}
\label{eq5.15}
\sup_{0 \leq t < \infty}\frac{|W_1(t)-t W_2(1)|}{(1+t)(\frac{t}{1+t})^\gamma}\overset{\cal L} =\sup_{ 0 \leq y \leq \frac{1}{{\cal D}^2}} |W(y)| \frac{(1+y-{\cal D}^2y)^\gamma}{y^\gamma}.
\end{equation}
By the asymptotic properties of a nonlinear regression, we have that the variance error estimator $\hat \sigma^2_m$ is strongly converging to $\sigma^2$, $|\hat \sigma_m - \sigma |=o_{\eP}(1)$.
The assertion \textit{(i)}  follows by the last relation together the relations (\ref{F1}), (\ref{LM}), (\ref{eq5.14})-(\ref{eq5.15}).\\
\textit{(ii)} The proof is similar of \textit{(i)}. We give its outline:
\[
\sup_{1 \leq k \leq T_m}\sum^{m+k}_{i=m+1} \hat \varepsilon_i/g(m,k,\gamma)  =\sigma \sup_{1 \leq k \leq T_m} [W_{1,m}(k)-{\cal D} \frac{k}{m} W_2(m)]/g(m,k,\gamma)(1+o_{\eP}(1))
 \]
 \[
\qquad \qquad \qquad \overset{{a.s.}} {\underset{m \rightarrow \infty}{\longrightarrow}} \sup_{0 \leq t \leq T} \frac{|W_1(t) -{\cal D}t W_2(1) |}{(1+t) \pth{\frac{t}{1+t}}^\gamma} \overset{\cal L} = \sup_{ 0 \leq y \leq \frac{T}{1+{\cal D}^2 T}} |W(y)| \frac{(1+y-{\cal D}^2y)^\gamma}{y^\gamma}.
 \]
\hspace*{\fill}$\blacksquare$\\

\noindent {\bf Proof of Theorem  \ref{theorem 2.2}}\\
We choose this particular  $k$: $\tilde k_m=k^0_m+m$. We will prove that for this $\tilde k_m$ we have $\lim_{m \rightarrow \infty}  \left| \sum^{m+\tilde k_m}_{i=m +1}\hat \varepsilon_i \right|/g(m, \tilde k_m,\gamma) =\infty$. Let us consider the partial sum of the residuals after the first $m$ observations
\begin{equation}
\label{et}
\sum^{m+\tilde k_m}_{i=m+1} \hat \varepsilon_i =\sum^{m+\tilde k_m}_{i=m+1} \varepsilon_i +\sum^{m+\tilde k_m}_{i=m+1} [f(\X_i;\ebo)-f(\X_i; \hat \eb_m)]+\sum^{m+\tilde k_m}_{i=m+k^0_m+1} [f(\X_i;\eb^0_m)-f(\X_i;\ebo)].
\end{equation}
Similar as for the Theorem \ref{theorem 2.1} we have, for the first two terms of the right-hand side of (\ref{et}),  
\begin{equation}
\label{ett}
\left|\sum^{m+\tilde k_m}_{i=m+1}  [\varepsilon_i+ f(\X_i;\ebo) - f(\X_i;\hat \eb_m)] \right|/g(m,\tilde k_m,\gamma)=O_{\eP}(1)
\end{equation} 
and for the last term of the right-hand side of (\ref{et})
\[
\sum^{m+\tilde k_m}_{i=m+k^0_m+1} [f(\X_i;\eb^0_m)-f(\X_i;\ebo)]=\sum^{m+\tilde k_m}_{i=m+k^0_m+1}[f(\X_i;\eb^0_m)-\eE[f(\X;\eb^0_m)]]
\]
\begin{equation}
\label{C1}
- \sum^{m+\tilde k_m}_{i=m+k^0_m+1} [f(\X_i;\ebo)-\eE[f(\X;\ebo)]]+(\tilde k_m-k^0_m) \pth{\eE[f(\X;\eb^0_m)]-\eE[f(\X;\ebo)]}.
\end{equation}
Since  $\eE[f(\X;\eb^0_m)] \neq \eE[f(\X;\ebo)]$, which implies, for the third  term of the right-hand side of (\ref{C1}) that
\begin{equation}
\label{C2}
\frac{(\tilde k_m-k^0_m) |\eE[f(\X;\eb^0_m)]-\eE[f(\X;\ebo)]|}{m^{1/2} g(m, \tilde k_m,\gamma)}=\frac{C  m}{m (1+\frac{\tilde k_m}{m})(\frac{\tilde k_m/m}{1+{\tilde k_m}/{m}})^\gamma}>C>0,
\end{equation} 
where $C$ is a  constant  not depending of $m$. For the last relation, we have used that for $x >1$ we have $\frac{1}{2} < \frac{x}{1+x} <1$, then $(\frac{\tilde k_m/m}{1+{\tilde k_m}/{m}})^\gamma \in (2^{-\gamma},1)$ and $(1+x)^{-1} \geq 1$. On the other hand, using assumption (A5)
\[
\sum^{m+\tilde k_m}_{i=m+k^0_m+1}[f(\X_i;\eb^0_m)-\eE[f(\X;\eb^0_m)]]=\sum^{m+\tilde k_m}_{i=1}[f(\X_i;\eb^0_m)-\eE[f(\X;\eb^0_m)]]- \sum^{m+k^0_m}_{i=1}[f(\X_i;\eb^0_m)-\eE[f(\X;\eb^0_m)]]
\]
is of order $ O_{\eP}(m+\tilde k_m)^{1/2}+O_{\eP}(m+k^0_m)^{1/2}= O_{\eP}(m+\tilde k_m)^{1/2}$. Moreover $ m^{-1}(1+\tilde k_m/m)^{-1}(m+\tilde k_m)^{1/2} \rightarrow 0$, as $m \rightarrow \infty$  and $(\frac{\tilde k_m/m}{1+{\tilde k_m}/{m}})^\gamma \in (2^{-\gamma},1)$. Thus, for the first term of the right-hand side of (\ref{et}) we have
\begin{equation}
\label{C3}
m^{-1/2} \sum^{m+\tilde k_m}_{i=m+k^0_m+1} [f(\X_i;\eb^0_m)-\eE[f(\X;\eb^0_m)]]/g(m,\tilde k_m,\gamma)=o_{\eP}(1).
\end{equation} 
Similarly, for the second term of the right-hand side of (\ref{C1})
\begin{equation}
\label{C4}
m^{-1/2} \sum^{m+\tilde k_m}_{i=m+k^0_m+1} [f(\X_i;\ebo)-\eE[f(\X;\ebo)]]/g(m,\tilde k_m,\gamma)=o_{\eP}(1).
\end{equation} 
Taking into account the relations (\ref{C1})-(\ref{C4}) we can get, for the third term of the right-hand side of (\ref{et})
\begin{equation}
\label{ettt}
\liminf_{m \rightarrow \infty} m^{-1/2}  \left|\sum^{m+\tilde k_m}_{i=m+k^0_m+1} [f(\X_i;\eb^0_m)-f(\X_i;\ebo)]\right| /g(m,\tilde k_m,\gamma)>0.
\end{equation}
The relations (\ref{et}), (\ref{ett}), (\ref{ettt}) imply that we found one $\tilde k_m$ such that $\lim_{m \rightarrow \infty} m^{-1/2} \left| \sum^{m+\tilde k_m}_{i=m +1}\hat \varepsilon_i \right|/g(m, \tilde k_m,\gamma) >0$. Thus $\lim_{m \rightarrow \infty}  \left| \sum^{m+\tilde k_m}_{i=m +1}\hat \varepsilon_i \right|/g(m, \tilde k_m,\gamma) =\infty$. Then
$
\lim_{m \rightarrow \infty} \sup_{1 \leq k < \infty} \left| \sum^{m+k}_{i=m +1}\hat \varepsilon_i \right|/g(m, k,\gamma) =\infty$. The theorem follows. 
\hspace*{\fill}$\blacksquare$\\

\noindent {\bf Proof of Proposition \ref{Proposition A}}\\
It is clear that
$|Z_{k,n}|-\mu_{k,n} \leq |Z_{k,n} -\mu_{k,n}| \leq \max_k |Z_{k,n} -\mu_{k,n}|$.
Then $ \eP[\max_k |Z_{k,n} -\mu_{k,n}| \geq \epsilon] \geq \eP[|Z_{k,n}| -|\mu_{k,n}|\geq \epsilon]=\eP[|Z_{k,n}| \geq \epsilon +|\mu_{k,n}|] \geq \eP[\max_k(|Z_{k,n}|) \geq \epsilon +|\mu_{k,n}|]\geq \eP[\max_k(|Z_{k,n}|) \geq 2\epsilon ]$. For the last inequality we have used: for all $\epsilon>0$ there exists a natural number $n_\epsilon$ such that for all $n \geq n_\epsilon$ we have $|\mu_{k,n}| \leq \epsilon$.
\hspace*{\fill}$\blacksquare$\\

\noindent {\bf Proof of Proposition \ref{Lemma 4}}\\
We denote by $e_2 \equiv -l/m \sum^m_{j=1}\varepsilon_{{\cal U}_{m,k}(j)}$ and we remind the notation $D_A \equiv \A^t\A$. Without loss of generality, we take $l \leq k$, the other cases are similar. Consider now  the following random variable,  for $i=m+1, \cdots, m+l$,
\[
G_{i,k} \equiv D^{-1}_A\frac{1}{m}\pth{\sum^m_{j=1} \ef^t(\X_j;\ebo) \varepsilon_{{\cal U}_{m,k}(j)} }\ef(\X_i;\ebo)-\frac{1}{m}\sum^m_{j=1}\varepsilon_{{\cal U}_{m,k}(j)}.
\]
Consequently,  $-I_2+e_2=\sum^{m+l}_{i=m+1}G_{i,k}$.  The conditional expectation of $-I_2+e_2$ is
\[
\eE^*_{k,m}[-I_2+e_2]=\frac{1}{m} \sum^m_{j=1}\ef^t(\X_j;\ebo)\frac{1}{m+k} \sum^{m+k}_{i=1} \varepsilon_i\B^{-1}_{m}\c_1(m,k,l)-\frac{l}{m} \sum^m_{j=1}\frac{1}{m+k} \sum^{m+k}_{i=1} \varepsilon_i
\]
\[
\qquad =\bar \varepsilon_{m+k}\acc{ D^{-1}_A \cro{\frac{m+l}{m} \sum^m_{j=1}\ef^t(\X_j;\ebo) \frac{1}{m+l}\sum^{m+l}_{i=1} \ef(\X_i;\ebo)-\frac{m}{m} \sum^m_{j=1}\ef^t(\X_j;\ebo) \frac{1}{m}\sum^{m}_{i=1} \ef(\X_i;\ebo)} -l }.\]
On the other hand, by assumption (A1) for all $\epsilon>0$, there exists $M_1>0$ such that $\eP[m^{1/2}|\bar \varepsilon_{m+k} | >M_1]<\epsilon$. Thus, taking also into account the assumption (A4) for $\ef(\X_i;\ebo)$,  we get
\[
\eE^*_{k,m} \cro{ \frac{-I_2+e_2}{g(m,l,\gamma)}}=\frac{o_{\eP}(m+l) \bar \varepsilon_{m+k}}{m^{1/2} \pth{1+\frac{l}{m}} \pth{\frac{l/m}{1+l/m}}^\gamma}= \frac{o_{\eP}(1+l/m)m^{1/2} \bar \varepsilon_{m+k}}{\pth{1+\frac{l}{m}} \pth{\frac{l/m}{1+l/m}}^\gamma}.
\]
Using the relation (\ref{Km}), the last relation is  $o_{\eP}(1)m^{1/2} \bar \varepsilon_{m+k}=o_{\eP}(1)O_{\eP}(1)=o_{\eP}(1)$, for all $ l,  k=1, \cdots, T_m$. 
Consequently 
\[
\eE^*_{k,m} \cro{ \frac{-\sum^{m+l}_{i=m+1}G_{i,k}}{g(m,l,\gamma)}}=o_{\eP}(1).
\]
for all $ l,  k=1, \cdots, T_m$. Then, we are in the conditions to apply the inequality (\ref{CC}) for the random variable $G_{i,k}$ and the sequence $(b_l \equiv g(m,l,\gamma))$. Hence, for any $\epsilon > 0$, there exists a natural number $m_\epsilon$ such that for $m \geq m_\epsilon$,
\begin{equation}
\label{eq32bis}
\eP^*_{k,m} \cro{\sup_{1 \leq l \leq T_m}  \frac{|\sum^{m+l}_{i=m+1} G_{i,k} |}{g(m,l,\gamma)} \geq 2 \epsilon} \leq \frac{1}{\epsilon^2} \sum^{T_m}_{l=1} \frac{\eE^*_{k,m}[G_{i,k}^2]}{g^2(m,l,\gamma)}.
\end{equation}
By elementary algebra, using the fact that for  $j \neq j'$, $\eE^*_{k,m} [ \varepsilon_{{\cal U}_{m,k}(j)}\cdot \varepsilon_{{\cal U}_{m,k}(j')} ]=\eE^*_{k,m} [ \varepsilon_{{\cal U}_{m,k}(j)} ] \cdot \eE^*_{k,m} [ \varepsilon_{{\cal U}_{m,k}(j')} ]=(m+k)^{-2} \pth{\sum^{m+k}_{a=1}\varepsilon_a}^2$,  yield
\begin{equation}
\label{sgl0}
\begin{array}{cl}
\eE^*_{k,m}[G_{i,k}^2]=&( \bar \varepsilon_{m+k})^2 \cro{m^{-1}\sum^m_{j=1} \pth{D^{-1}_A\ef^t(\X_j;\ebo) \ef(\X_i;\ebo)  -1 } }^2
\\
 &+\pth{\overline{\varepsilon^2}_{m+k}-( \bar \varepsilon_{m+k})^2} \cro{ m^{-2}\sum^m_{j=1}  \pth{ D^{-1}_A\ef^t(\X_j;\ebo) \ef(\X_i;\ebo) -1}^2 }.
\end{array}
\end{equation}
By the relation (30) of Hu\v{s}kov\'a and Kirch (2012), we get, for a constant $C_1 >0$: $g(m,l,\gamma) \geq  C_1  ( m^{1/2-\gamma}l^\gamma \e1_{l \leq m} + m^{-1/2}l \e1_{l>m} ) $. 
Thus, for  $l>m$ we have
$
g^{-2}(m,l,\gamma) < mC_1^{-2}l^{-2}< C_1^{-2} m^{-1} \rightarrow 0$, as $m \rightarrow \infty$ 
and for $1 \leq l \leq m$, $g^{-2}(m,l,\gamma) \leq C^{-2}_1 m^{-1+2 \gamma} l^{-2 \gamma} \leq C_1^{-2} m^{-1+2 \gamma} \rightarrow 0$,  as $m \rightarrow \infty$. 
 Under the assumptions (A4) and (A6) we have the following inequalities with a probability close to 1
\[
\frac{1}{m^2}\sum^m_{l=1} \cro{\sum^m_{j=1}  \pth{D^{-1}_A\ef^t(\X_j;\ebo) \ef(\X_l;\ebo) -1 } }^2\frac{1}{g^2(m,l,\gamma)} \leq C \sum^m_{l=1}\frac{1}{g^2(m,l,\gamma)} \leq C .
\]
On the other hand, by the Cauchy-Schwarz inequality we readily have with a probability 1
\[
 \frac{1}{m^2} \sum^{T_m}_{l=m+1} \cro{\sum^m_{j=1} \pth{\ef^t(\X_j;\ebo) \ef(\X_l;\ebo) D^{-1}_A -1 } }^2\frac{1}{g^2(m,l,\gamma)}\]
 \[ \qquad \qquad   \leq  \acc{\frac{1}{m^4} \sum^{T_m}_{l=m+1} \cro{\sum^m_{j=1} \pth{\ef^t(\X_j;\ebo) \ef(\X_l;\ebo) D^{-1}_A -1 } }^4}^{1/2} \acc{\sum^{T_m}_{l=m+1}\frac{1}{g^4(m,l,\gamma)} }^{1/2} .
\]
\[
\leq \acc{C (T_m-m)\sum^{T_m}_{l=m+1}\frac{1}{g^4(m,l,\gamma)}  }^{1/2} \leq \acc{C (T_m-m) m^2 \sum^{T_m}_{l=m+1}\frac{1}{l^4}  }^{1/2} =C \acc{m^2(T_m-m) \pth{\frac{1}{m^4}- \frac{1}{T^4_m} }}^{1/2}=C .
\] 
were used that $ m^{-1}\sum^m_{j=1}\ef^t(\X_j;\ebo) $, $ m^{-1}\sum^m_{j=1}\ef^t(\X_j;\ebo) \ef^t(\X_j;\ebo)$,  $ T_m^{-1}\sum^{T_m}_{j=1}\ef^t(\X_j;\ebo) $ are converging by assumption (A4). Hence
\begin{equation}
\label{sgl1}
\frac{1}{m^2}\sum^{T_m}_{l=1} \cro{\sum^m_{j=1} \pth{D^{-1}_A\ef^t(\X_j;\ebo) \ef(\X_l;\ebo)  -1 } }^2\frac{1}{g^2(m,l,\gamma)}=O_{\eP}(1) .
\end{equation}
For the second term of the right-hand side  of the  relation (\ref{sgl0}) we have: $\ef^t(\X_j;\ebo)\ef(\X_l;\ebo) \ef^t(\X_j;\ebo)\ef(\X_l;\ebo)$\\ $=trace(\ef^t(\X_l;\ebo)\ef(\X_j;\ebo) \ef^t(\X_j;\ebo)\ef(\X_l;\ebo))=\ef^t(\X_l;\ebo)\ef(\X_j;\ebo) \ef^t(\X_j;\ebo)\ef(\X_l;\ebo)$. Consequently, since for  $l \leq m$, $g^{-2}(m,l,\gamma) \leq C m^{-1+2\gamma}$, we have with a probability close to 1 that $m^{-2}\sum^m_{l=1}\sum^m_{j=1}  \{ \ef^t(\X_j;\ebo) \ef(\X_l;\ebo) D^{-1}_A-1\}^2 g^{-2}(m,l,\gamma)$ is less than or equal to
$2m^{-1}\sum^m_{l=1} \cro{\ef^t(\X_l;\ebo) \\B_m \ef(\X_l;\ebo) D^{-2}_A+1 } g^{-2}(m,l,\gamma)$
$\leq Cm^{-1+2\gamma} m^{-1}$ \\ $\cdot \sum^m_{l=1} [\ef^t(\X_l;\ebo)\B_m  \ef(\X_l;\ebo) D^{-2}_A+1 ]\leq C m^{-1+2\gamma} m^{-1} \sum^m_{l=1} [\| D^{-2}_A\ef^t(\X_l;\ebo)\|^2_2 \|\B_m\|_2 +1]$ 
$ =C m^{-1+2 \gamma}\rightarrow 0$, for $m \rightarrow \infty$. For the last relation we used the assumption (A4). \\
We have in the other hand
$m^{-2}\sum^{T_m}_{l=m+1}\sum^m_{j=1}  \pth{ \ef^t(\X_j;\ebo) \ef(\X_l;\ebo) D^{-1}_A-1}^2 g^{-2}(m,l,\gamma)$ is less than or equal to
$
2m^{-1}\sum^{T_m}_{l=m+1} g^{-2}(m,l,\gamma) m^{-1} \sum^m_{j=1} \cro{D^{-2}_A \ef^t(\X_l;\ebo)\ef(\X_j;\ebo) \ef^t(\X_j;\ebo)\ef(\X_l;\ebo)+1}
$\\
$
=2m^{-1}\sum^{T_m}_{l=m+1} g^{-2}(m,l,\gamma) [D^{-2}_A \ef^t(\X_l;\ebo) \B_m \ef(\X_l;\ebo)+1  ]
$
and by the  Cauchy-Schwarz inequality\\
$
\leq 2m^{-1} \acc{\sum^{T_m}_{l=m+1} \cro{\| D^{-2}_A\ef(\X_l;\ebo)\|^2_2 \|\B_m\|_2 +1}^2 }^{1/2} \acc{\sum^{T_m}_{l=m+1} g^{-4}(m,l,\gamma)}^{1/2}
$
\begin{equation}
\label{sgl2}
\leq Cm^{-1} \acc{(T_m-m) \sum^{T_m}_{l=m+1} g^{-4}(m,l,\gamma) }^{1/2}\leq Cm^{-1}\rightarrow 0, \quad \textrm{for } m \rightarrow \infty,
\end{equation}
uniformly in $k$. Since $ \bar \varepsilon_{m+k} \overset{{a.s}} {\underset{m \rightarrow \infty}{\longrightarrow}} 0$ and $ \overline{\varepsilon}^2_{m+k}- (\bar \varepsilon_{m+k})^2 =\hat \sigma^2_{m,k} \overset{{a.s}} {\underset{m \rightarrow \infty}{\longrightarrow}}  \sigma^2$, and using the results (\ref{sgl0}), (\ref{sgl1}), (\ref{sgl2}), we obtain by  (\ref{eq32bis}) that
$
\sup_{1 \leq k < \infty} \eP^*_{k,m} \cro{ \sup_{1 \leq l \leq T_m} \left|\sum^{m+l}_{i=m+1} G_i \right|/g(m,l,\gamma) \geq \epsilon}\overset{{\eP}} {\underset{m \rightarrow \infty}{\longrightarrow}} 0$.
\hspace*{\fill}$\blacksquare$\\

\noindent {\bf Proof of Proposition  \ref{lemma 5}}\\
Using the Proposition \ref{Lemma 4} and the   Lemmas \ref{Lemma I3I6}, \ref{Lemma I4I7I5I8},  the proof is similar to that of the  Lemma 5 of Hu\v{s}kov\'a and Kirch (2012).
\hspace*{\fill}$\blacksquare$\\

\noindent {\bf Proof of Theorem  \ref{theorem 1HK}}\\
Using the Proposition \ref{lemma 5}, the Theorems \ref{theorem 2.1}, \ref{theorem 2.2}  and Lemma \ref{lemma 7}, the proof is similar to that of the  Theorem 1 of Horv\'ath and Kirch (2012).
\hspace*{\fill}$\blacksquare$\\

\section{Appendix}
In this section useful  Lemmas to prove the main results of Sections 3 and 4 are given. 
 We recall the notations: $\A \equiv \eE[ \ef(\X;\ebo)]$,  $\B\equiv \eE[\ef(\X;\ebo) \ef^t(\X;\ebo)]$, $A_i\equiv \ef^t(\X_i;\ebo) \B^{-1} \A$.

\subsection{Lemmas for Section 3} 
 
\begin{lemma}
\label{Lemma 5.2} Suppose that assumptions  (A1)-(A4) hold. 
Under the hypothesis $H_0$ we have, as $m \rightarrow \infty$,
\[
\sup_{1 \leq k <\infty}\left|\sum^{m+k}_{i=m+1}\hat \varepsilon_i-\pth{\sum^{m+k}_{i=m+1} \varepsilon_i - \frac{k}{m} \sum^m_{i=1} A_i \varepsilon_i} \right|/g(m,k,\gamma)=o_{\eP}(1).
\]

\end{lemma}
\noindent {\bf Proof of Lemma \ref{Lemma 5.2}}\\
Under the hypothesis $H_0$, $\sum^{m+k}_{i=m+1}\hat \varepsilon_i=\sum^{m+k}_{i=m+1} \varepsilon_i - \sum^{m+k}_{i=m+1} [f(\X_i;\hat \eb_m)-f(\X_i;\ebo)]$. Then, by a Taylor expansion of $f(\X_i;\eb)$ in a neighborhood of $\ebo$
\begin{equation}
\label{eq10}
\sum^{m+k}_{i=m+1}\hat \varepsilon_i- \pth{\sum^{m+k}_{i=m+1} \varepsilon_i -\frac{k}{m}\sum^m_{i=1} A_i \varepsilon_i}
=\frac{k}{m}\sum^m_{i=1} A_i \varepsilon_i- \sum^{m+k}_{i=m+1} [f(\X_i;\hat\eb_m)-f(\X_i;\ebo)]
= \frac{k}{m}\sum^m_{i=1} A_i \varepsilon_i -(\hat \eb_m-\ebo)^t\sum^{m+k}_{i=m+1} \ef(\X_i; \tilde \eb_m),
\end{equation}
with $\tilde \eb_m=\ebo+\eth (\hat \eb_m-\ebo)$, $\eth \in [0,1]^q$.
We know that the LS estimator $\hat \eb_m$ of parameter in a linear model is $\sqrt{m}$-consistent (see Seber and Wild, 2003) $\hat \eb_m-\ebo=O_{\eP}(m^{-1/2})$. On the other hand, by the triangle inequality
\begin{equation}
\label{point}
\left\| \sum^{m+k}_{i=m+1}  \ef(\X_i;\tilde \eb_m)-k \eE[\ef(\X;\ebo)] \right\|_1  \leq  \left\| \sum^{m+k}_{i=1}\pth{ \ef(\X_i;\tilde \eb_m)- \eE[\ef(\X;\ebo)]}\right\|_1 \\
 +\left\|\sum^m_{i=1}\pth{ \ef(\X_i;\tilde \eb_m)- \eE[\ef(\X;\ebo)]} \right\|_1
\end{equation}
Generally, for any $n \geq m$ and  $\eb$ in a $m^{-1/2}$-neighborhood of  $\ebo$, by assumption (A4) for $\ef(\X_i;\ebo)$, using the law of iterated logarithm we have that for all $\epsilon>0$, there exists a $M_2 >0$ such that $\eP[n^{-1/2} \sum^n_{i=1} \| \ef(\X_i;\ebo)-\A\|_1 \geq M_2] \leq \epsilon$. Together with assumption (A2), it holds that, for all $\eb$ in a $m^{-1/2}$-neighborhood of $\ebo$ 
\[
\left\| \sum^n_{i=1}\cro{\ef(\X_i;\eb)-\eE[\ef(\X;\ebo)]}\right\|_1 \leq \left\| \sum^n_{i=1}\cro{\ef(\X_i;\eb)- \ef(\X_i;\ebo)}\right\|_1+\left\|\sum^n_{i=1}\cro{\ef(\X_i;\ebo) -\eE[\ef(\X;\ebo)]}\right\|_1
\]
$= \|\eb-\ebo\|_1O_{\eP}(n)+O_{\eP}(n^{1/2})=O_{\eP}(n \cdot m^{-1/2})+O_{\eP}(n^{1/2})$ uniformly in $\eb$.
Thus, the right-hand side of  (\ref{point}) becomes $(m+k)m^{-1/2}+(m+k)^{1/2}+m^{1/2}+m^{1/2}=(m+k)m^{-1/2}+(m+k)^{1/2}+2m^{1/2}$.  Hence, for the last term of (\ref{eq10}) we have
\[
\sup_{1 \leq k < \infty}\frac{\left|(\hat \eb_m-\ebo)^t\sum^{m+k}_{i=m+1}\pth{\ef(\X_i;\tilde \eb_m)- \eE[\ef(\X;\ebo)] } \right|}{g(m,k,\gamma)}=O_{\eP}(m^{-1/2})\sup_{1 \leq k < \infty} \frac{m^{1/2}+(m+k)m^{-1/2}+(m+k)^{1/2}}{m^{1/2}\pth{1+\frac{k}{m}} \pth{\frac{k}{m+k}}^\gamma}
\]
\[
= C  \sup_{1 \leq k < \infty} \frac{1+\pth{1+\frac{k}{m}} +\pth{1+\frac{k}{m}}^{1/2}}{m^{1/2}\pth{1+\frac{k}{m}} \pth{\frac{k}{m+k}}^\gamma} \leq C \sup_{1 \leq k < \infty} \frac{2+\frac{k}{m}}{m^{1/2}\pth{1+\frac{k}{m}} \pth{\frac{k}{m}\cdot \frac{1}{1+\frac{k}{m}}}^\gamma}.
\]
Since, for $k \leq m$, $\frac{k}{k+m} > \frac{1}{2m}$ and for $x>0 $ we have $(1+x)^{-1} <1$, we can write
\[
\sup_{1 \leq k \leq m} \frac{2+\frac{k}{m}}{m^{1/2}\pth{1+\frac{k}{m}} \pth{\frac{k}{m}\cdot \frac{1}{1+\frac{k}{m}}}^\gamma} \leq \frac{3}{m^{1/2}\pth{\frac{1}{2m}}^\gamma}=3 \cdot 2^{\gamma-1} m^{-1/2+\gamma}=o(1).
\]
For all $x \geq 1$, we have that $\pth{\frac{x+1}{x}}^\gamma \leq 2^\gamma$ and $(1+x)^{-1} \leq 2^{-1}$. Then
\[
 \sup_{m \leq k < \infty} \frac{2+\frac{k}{m}}{m^{1/2}\pth{1+\frac{k}{m}} \pth{\frac{k}{m}\cdot \frac{1}{1+\frac{k}{m}}}^\gamma} \leq \frac{1}{m^{1/2}}\pth{\frac{2^\gamma}{2}+2^\gamma }=o(1).
\]
Hence, (\ref{eq10}) becomes
\begin{equation}
\label{eq11}
\sum^{m+k}_{i=m+1}\hat \varepsilon_i- \pth{\sum^{m+k}_{i=m+1}\varepsilon_i-\frac{k}{m}\sum^m_{i=1}A_i\varepsilon_i}=\frac{k}{m}\sum^m_{i=1}A_i\varepsilon_i+k(\hat \eb_m-\ebo)^t \eE[\ef(\X;\ebo)](1+o_{\eP}(1)) .
\end{equation}
On the other hand, $\hat \eb_m$ is the least squares estimator of $\ebo$, calculated for $i=1, \cdots, m$,
\[
\begin{array}{lll}
0 & = & \sum^m_{i=1} \ef(\X_i;\hat \eb_m)[Y_i-f(\X_i;\hat \eb_m)]= \sum^m_{i=1} \ef(\X_i;\hat \eb_m) [\varepsilon_i - (\hat \eb_m-\ebo)^t \ef(\X_i;\tilde \eb_m)]\\
 & = & \sum^m_{i=1} \ef(\X_i;\ebo) \varepsilon_i + \sum^m_{i=1} \eff(\X_i; \tilde \eb_m) (\hat \eb_m-\ebo)^t\varepsilon_i -\sum^m_{i=1} \ef(\X_i;\ebo) \ef^t(\X_i;\ebo)(\hat \eb_m-\ebo)\\
 & & - 1/2(\hat \eb_m-\ebo)^t \sum^m_{i=1}\eff(\X_i;\tilde \eb_m) \ef(\X_i;\ebo) (\hat \eb_m-\ebo). 
% & = & \eE[\ef(\X;\ebo)]\sum^m_{i=1} \varepsilon_i+\sum^m_{i=1}\cro{\ef(\X_i;\hat \eb_m)-\ef(\X_i;\ebo)} \varepsilon_i +\sum^m_{i=1} \cro{\ef(\X_i;\ebo)-\eE[\ef(\X;\ebo)]} \varepsilon_i\\
%  & & -m(\hat \eb_m-\ebo)^t \eE[\ef(\X;\ebo) \ef^t(\X;\ebo)]-(\hat \eb_m-\ebo)\sum^m_{i=1} \cro{\ef(\X_i;\hat \eb_m) \ef^t(\X_i;\tilde \eb_m)- \eE[\ef(\X;\ebo) \ef^t(\X;\ebo)]}\\
 % & =& \sum^m_{i=1} \ef(\X_i;\ebo)\varepsilon_i -m(\hat \eb_m-\ebo) \eE[\ef(\X;\ebo) \ef^t(\X;\ebo)] {1+o_{\eP}(1)}.
\end{array}
\]
Using the assumptions  (A1), (A2) and the  Cauchy-Schwarz inequality, we obtain
\begin{equation}
\label{etoile}
0=\cro{\sum^m_{i=1} \ef(\X_i;\ebo) \varepsilon_i  -\sum^m_{i=1} \ef(\X_i;\ebo) \ef^t(\X_i;\ebo)(\hat \eb_m-\ebo) }(1+o_{\eP}(1)).
\end{equation}
Then, by relation (\ref{etoile}) below
\begin{equation}
\label{estbeta}
\hat \eb_m-\ebo=\B^{-1}_m\pth{\frac{1}{m} \sum^m_{i=1} \ef(\X_i;\ebo) \varepsilon_i} (1+o_{\eP}(1))=\B^{-1}\pth{\frac{1}{m} \sum^m_{i=1} \ef(\X_i;\ebo) \varepsilon_i} (1+o_{\eP}(1)),
\end{equation} 
again too
\[
k(\hat \eb_m-\ebo)^t \eE[\ef(\X_i;\ebo)] =\pth{\frac{k}{m} \sum^m_{i=1} \ef^t(\X_i;\ebo) \varepsilon_i} \B^{-1} \eE[\ef(\X;\ebo)](1+o_{\eP}(1)) =\frac{k}{m} \sum^m_{i=1} A_i \varepsilon_i (1+o_{\eP}(1))
\] and  we  replace next in (\ref{eq11}). To complete the proof, we must prove that for  (\ref{etoile}) that 
 $\sup_{1 \leq k <\infty} k/(m g(m,k,\gamma)) o_{\eP}(1)=o_{\eP}(1)$. Using (A1)-(A4) and  the fact that $ \hat \eb_m - \ebo=O_{\eP}(m^{-1/2})$  we deduce that
\[
\sup_{1 \leq k <\infty} \frac{\| (\hat \eb_m-\ebo)^t \frac{k}{m} \sum^m_{i=1} \eff(\X_i;\tilde \eb_m) \varepsilon_i \|_1}{m^{1/2} \pth{1+\frac{k}{m}} \pth{\frac{k}{m+k}}^\gamma} \leq K_m \frac{1}{m} \| \sum^m_{i=1} \eff(\X_i;\tilde \eb_m) \varepsilon_i \|_1 =o_{\eP}(1),
\]
with $K_m$ given by (\ref{Km}). Similarly
\[
\sup_{1 \leq k <\infty} \frac{\| (\hat \eb_m-\ebo)^t \frac{k}{m} \sum^m_{i=1} \{\ef(\X_i; \ebo) \ef^t(\X_i;\ebo)- \eE[\ef(\X;\ebo) \ef^t(\X;\ebo)]\}\|_1}{g(m,k,\gamma)}=K_m o_{\eP}(1)=o_{\eP}(1).
\]
Using (\ref{Km}), with the Cauchy-Schwarz inequality for matrix, we have
\[
\sup_{1 \leq k <\infty} \frac{\| (\hat \eb_m-\ebo)^t \sum^m_{i=1}\eff(\X_i;\tilde \eb_m) \ef(\X_i;\ebo) (\hat \eb_m-\ebo)^t \|_2}{g(m,k,\gamma)}\]
\[\qquad \qquad \leq \sup_{1 \leq k <\infty} K_m \|\frac{k}{m} \sum^m_{i=1}\eff(\X_i;\tilde \eb_m) \ef(\X_i;\ebo) \|_1 \cdot \| \hat \eb_m-\ebo \|_1 \cdot O_{\eP}(\| \hat \eb_m-\ebo \|_1)=o_{\eP}(1).
\]
\hspace*{\fill}$\blacksquare$\\

\begin{lemma}
\label{Lemma 5.3} Suppose that assumptions  (A1)-(A3) hold. 
Under the hypothesis $H_0$,  there exists two independent Wiener processes $\{W_{1,m}(t), 0 \leq t < \infty \}$ and $\{W_{2,m}(t), 0 \leq t < \infty \}$ such that, for $m \rightarrow \infty$,
\[ 
\sup_{1 \leq k <\infty} \left| \pth{\sum^{m+k}_{i=m+1} \varepsilon_i- \frac{k}{m} \sum^m_{i=1}A_i \varepsilon_i} - \pth{ \sigma W_{1,m} (k) -\frac{k}{m} \sigma {\cal D} W_{2,m}(m)} \right| /g(m,k,\gamma)=o_{\eP}(1).
\]
\end{lemma}
\noindent {\bf Proof of Lemma \ref{Lemma 5.3}}\\
The random variables $\{\sum^{m+k}_{i=m+1} \varepsilon_i, 1 \leq k < \infty\}$ and $\{\sum^m_{i=1}A_i \varepsilon_i \}$ are  independent. It is obvious that, since  $\X_i$ is  independent of $\varepsilon_i$, we have $\eE[A_i \varepsilon_i]=0$, $Var[A_i \varepsilon_i]=\eE[A_i^2]\eE[\varepsilon^2_i]$. On the other hand, $\eE[\varepsilon^2_i]=\sigma^2$  and $\eE[A_i^2]=\A^t\B^{-1}\A$. By an argument similar to the one used in  Horv\'ath et al.(2004), Lemma 5.3., we obtain
$\sup_{1 \leq k <\infty} \left|\sum^{m+k}_{i=m+1} \varepsilon_i- \sigma W_{1,m}(k) \right|/k^{1/\gamma}=O_{\eP}(1) $ and   $\sum^m_{i=1}A_i \varepsilon_i -\sigma {\cal D} W_{2,m} (m) =o_{\eP}(m^{1/\nu}) $, as $m \rightarrow \infty$, $\nu >2$. The rest of proof is similar to that of the Lemme 5.3. of Horv\'ath et al.(2004).
\hspace*{\fill}$\blacksquare$\\

\subsection{Lemmas for Section 4}
We recall that (see the decomposition of $g \tilde \Gamma$ given in Section 4): $g(m,l,\gamma) \tilde \Gamma(m,l,\gamma)( \varepsilon^*_{m,k}(1), \cdots, \varepsilon^*_{m,k}(m+l) )=I_1+I_2+{\cal R}_m$. More precisely, ${\cal R}_m$ have the decomposition: ${\cal R}_m \equiv \sum^8_{j=3}I_j$, with
\[
I_3 \equiv \sum^{m+l}_{i=m+1} \e1_{k \leq k^0_m} \e1_{{\cal U}_{m,k}(i) < m+k^0_m}\cro{f(\X_{{\cal U}_{m,k}(i)};\ebo)-f(\X_{{\cal U}_{m,k}(i)};\hat \eb_{m+k}) },
\]
\[
I_4 \equiv \sum^{m+l}_{i=m+1} \e1_{k > k^0_m} \e1_{{\cal U}_{m,k}(i) < m+k^0_m}\cro{f(\X_{{\cal U}_{m,k}(i)};\ebo)-f(\X_{{\cal U}_{m,k}(i)};\hat \eb_{m+k}) },
\]
\[
I_5 \equiv \sum^{m+l}_{i=m+1} \e1_{k > k^0_m} \e1_{{\cal U}_{m,k}(i) > m+k^0_m}\cro{f(\X_{{\cal U}_{m,k}(i)};\eb^0_m)-f(\X_{{\cal U}_{m,k}(i)};\hat \eb_{m+k}) },
\]
\[
I_6 \equiv-\pth{\frac{1}{m} \sum^m_{j=1} \ef^t(\X_j;\ebo)\e1_{k \leq k^0_m} \e1_{{\cal U}_{m,k}(j) < m+k^0_m}\cro{f(\X_{{\cal U}_{m,k}(j)};\ebo)-f(\X_{{\cal U}_{m,k}(j)};\hat \eb_{m+k}) }}\B^{-1}_{m} \c_1(m,k,l) ,
\]
\[
I_7 \equiv- \pth{\frac{1}{m} \sum^m_{j=1} \ef^t(\X_j;\ebo)\e1_{k > k^0_m} \e1_{{\cal U}_{m,k}(j) < m+k^0_m}\cro{f(\X_{{\cal U}_{m,k}(j)};\ebo)-f(\X_{{\cal U}_{m,k}(j)};\hat \eb_{m+k}) }}\B^{-1}_{m} \c_1(m,k,l),
\]
\[
I_8 \equiv- \pth{\frac{1}{m} \sum^m_{j=1} \ef^t(\X_j;\ebo)\e1_{k > k^0_m} \e1_{{\cal U}_{m,k}(j) > m+k^0_m}\cro{f(\X_{{\cal U}_{m,k}(j)};\eb^0_m)-f(\X_{{\cal U}_{m,k}(j)};\hat \eb_{m+k}) }}\B^{-1}_{m} \c_1(m,k,l).
\]
Under the hypothesis $H_0$, $I_4=I_5=I_7=I_8=0$. Let us consider now
\[
\tilde I_3 \equiv I_3 -\frac{l}{m+k} \sum^{m+k}_{j=1} [f(\X_j;\ebo)- f(\X_j;\hat \eb_{m+k})],
\quad
\tilde I_6 \equiv I_6+ \frac{l}{m+k} \sum^{m+k}_{j=1} [f(\X_j;\ebo)- f(\X_j;\hat \eb_{m+k})].
\]

\begin{lemma}
\label{Lemma I3I6}
Under the assumptions (A1)-(A4), for all $\epsilon >0$, we have
\begin{equation}
\label{eq_I3}
\eP^*_{m,k} \cro{\max_{1 \leq l \leq T_m} \frac{| \tilde I_3|}{g(m,l,\gamma)} \geq \epsilon} \rightarrow 0 \qquad \textrm{in probability, uniformly in  }k, \quad \textrm{as } m \rightarrow \infty, 
\end{equation}
\begin{equation}
\label{eq_I6}
\sup_{1 \leq k < \infty} \eP^*_{m,k} \cro{\max_{1 \leq l \leq T_m} \frac{|\tilde I_6|}{g(m,l,\gamma)} \geq \epsilon} \rightarrow 0, \qquad \textrm{in probability, as } m\rightarrow \infty, 
\end{equation}
whether under $H_0$ or $H_1$.
\end{lemma}
\noindent {\bf Proof of Lemma \ref{Lemma I3I6}}
Let us consider the random  variable, for $i=m+1, \cdots, m+l$,
\begin{equation}
\label{zz}
{\cal Z}_{k,m}(i)=f(\X_{{\cal U}_{m,k}(i)};\ebo)-f(\X_{{\cal U}_{m,k}(i)};\hat \eb_{m+k}) -\frac{1}{m+k} \sum^{m+k}_{j=1}[f(\X_j;\ebo)-f(\X_j;\hat \eb_{m+k})].
\end{equation}
Then $\eE^*_{m,k}[{\cal Z}_{k,m}(i)]=(m+k)^{-1}\sum^{m+k}_{j=1}  [f(\X_j;\ebo)- f(\X_j;\hat \eb_{m+k})] -(m+k)^{-1} \sum^{m+k}_{j=1}[f(\X_j;\ebo)-f(\X_j;\hat \eb_{m+k})]=0$. Hence, since $\tilde I_3=\sum^{m+l}_{i=m+1} {\cal Z}_{k,m}(i)$, we have $\eE^*_{m,k}[\tilde I_3]=0$.\\
For $\tilde I_6$ let us  consider $l \leq k$, for two other cases the arguments are like. Since $\tilde I_6$ is  scalar, using (A4), from an equality to an other we apply the  \textit{trace} operator, the conditional expectation $\eE^*_{m,k}[-I_6]$ is equal to 
\[
 \frac{1}{m+k} \sum^{m+k}_{i=1} [f(\X_i;\ebo)-f(\X_i;\hat \eb_{m+k})]\pth{\frac{1}{m} \sum^m_{j=1}\ef(\X_j;\ebo)}\B^{-1}_{m+k}\B_{m+k}  D^{-1}_A \pth{\sum^{m+l}_{i=m+1} \ef^t(\X_i;\ebo)} \qquad
\]
\[
\qquad \qquad =\frac{D^{-1}_A}{m+k} \sum^{m+k}_{i=1} [f(\X_i;\ebo)-f(\X_i;\hat \eb_{m+k})]\cro{\frac{m+l}{m+l}\sum^{m+l}_{i=1} \ef^t(\X_i;\ebo)-\frac{m}{m} \sum^m_{i=1}\ef^t(\X_i;\ebo)} \pth{\frac{1}{m} \sum^m_{i=1}\ef(\X_i;\ebo)} 
\]
\[
\qquad \;\;=[(m+l)\A^t-m\A^t]D^{-1}_A\A(1+o_{\eP}(1))\frac{1}{m+k} \sum^{m+k}_{i=1} [f(\X_i;\ebo)-f(\X_i;\hat \eb_{m+k})]
\]
\[
\qquad \;\; =l(1+o_{\eP}(1))\frac{1}{m+k} \sum^{m+k}_{i=1} [f(\X_i;\ebo)-f(\X_i;\hat \eb_{m+k})].
\]
Hence, 
$
\eE^*_{m,k}[\tilde I_6]= o_{\eP}(1)\frac{l}{m+k} \sum^{m+k}_{i=1} [f(\X_i;\ebo)-f(\X_i;\hat \eb_{m+k})]=l(m+k)^{-1/2}o_{\eP}(1)=lm^{-1/2}o_{\eP}(1)$.
On the other hand, the conditional variance of ${\cal Z}_{k,m}(i)$ is
\[
Var^*_{m,k}[{\cal Z}_{k,m}(i)] = \eE^*[ f(\X_{{\cal U}_{m,k}(i)};\eb^0_m)-f(\X_{{\cal U}_{m,k}(i)};\hat \eb_{m+k}) ]^2
=\frac{1}{m+k} \sum^{m+k}_{i=1} [f(\X_i;\ebo)-f(\X_i;\hat \eb_{m+k})]^2
\]
\[
\qquad \qquad \qquad \qquad \leq \frac{\|\hat \eb_{m+k}-\ebo\|_2^2}{m+k} \cro{\sum^{m+k}_{i=1} \|\ef(\X_i;\ebo)\|_2^2}(1+o_{\eP}(1))=C(m+k)^{-1} (1+o_{\eP}(1)),
\]
which is  $o_{\eP}(1)$ uniformly in  $k$. We apply the  H\'ajek-R\'enyi inequality for  ${\cal Z}_{k,m}(i)$: for all $\epsilon>0$,
\[
\eP^*_{m,k} \cro{\max_{1 \leq l \leq T_m} \frac{1}{g(m,l,\gamma)} \left|\sum^{m+l}_{i=m+1} {\cal Z}_{k,m}(i)\right| \geq \epsilon} \leq \frac{1}{\epsilon^2} \sum^{T_m}_{l=1} \frac{\eE^*_{m,k}[{\cal Z}^2_{k,m}(l) ]}{g^2(m,l,\gamma)}
\leq \frac{1}{\epsilon^2}\frac{C}{m+k}\sum^{T_m}_{l=1} \frac{1}{g^2(m,l,\gamma)}=\frac{C}{\epsilon^2(m+k)},
\]
and the relation (\ref{eq_I3}) follows.\\
For $\tilde I_6$, we can write, $\tilde I_6=-  \sum^m_{j=1} \ef^t(\X_j;\ebo){\cal Z}_{k,m}(j) \B^{-1}_{m} \c_1(m,k,l)$. Then, since
\[
\frac{\eE^*_{m,k}[\tilde I_6]}{g(m,l,\gamma)}=\frac{l/m}{\pth{1+l/m} \pth{\frac{l}{l+m}}^\gamma}o_{\eP}(1)=\pth{\frac{l/m}{1+l/m}}^{1-\gamma}  o_{\eP}(1) \quad \textrm{and } \pth{\frac{l/m}{1+l/m}}^{1-\gamma} \quad \textrm{is bounded by the  relation (\ref{Km}),}
\]
 combined with $\sup_{1 \leq k < \infty} \sup_{1 \leq l <\leq T_m }$ $\| \c_1(m,k,l)\|_2/m \leq C$ and since the random trial of  bootstrap are  independently, then ${\cal Z}_{k,m}(j)$ are also independently, we have
 $Var^*_{m,k}[\tilde I_6]= m^{-1} Var^*_{m,k}[{\cal Z}_{k,m}(j)]\c^t_1(m,k,l) \B^{-1}_{m+k} \c_1(m,k,l)=o_{\eP}(1)$.  
 Hence, by the   Bienaym\'e-Tchebychev  inequality,  Proposition \ref{Proposition A} and inequality (\ref{CC}),  the relation  (\ref{eq_I6}) follows.
\hspace*{\fill}$\blacksquare$\\

Let be (the expressions of $I_4, I_5, I_7, I_8$ are given before the Lemma \ref{Lemma I3I6}): 
\[\tilde I_4\equiv I_4 -\e1_{k >k^0_m} \frac{l}{m+k} \sum^{m+k^0_m}_{j=1} [ f(\X_j;\ebo)-f(\X_j;\hat \eb_{m+k})], \qquad \tilde I_5=I_5-\e1_{k >k^0_m} \frac{l}{m+k} \sum^{m+k}_{j=m+k^0_m+1} [ f(\X_j;\eb^0_m)-f(\X_j;\hat \eb_{m+k})],
\]
 \[\tilde I_7\equiv I_7 +\e1_{k >k^0_m} \frac{l}{m+k} \sum^{m+k^0_m}_{j=1} [ f(\X_j;\ebo)-f(\X_j;\hat \eb_{m+k})],
\qquad \tilde I_8=I_8+\e1_{k >k^0_m} \frac{l}{m+k} \sum^{m+k}_{j=m+k^0_m+1} [ f(\X_j;\eb^0_m)-f(\X_j;\hat \eb_{m+k})].
\]

\begin{lemma}
\label{Lemma I4I7I5I8}
Under the assumptions (A1)-(A5),  if the hypothesis $H_1$ is true, we have that for all $\epsilon >0$ there  exists $M>0$ such that we have in probability
\begin{equation}
\label{ineg_I4}
\sup_{1 \leq k < \infty} \eP^*_{m,k} \cro{\max_{1 \leq l \leq T_m} \frac{ |\tilde I_j|}{g(m,l,\gamma)} \geq M}  \leq \epsilon+o_{\eP}(1), \qquad \textrm{ for } j \in \{ 4,5,7,8\}.
\end{equation}
\end{lemma}

\noindent {\bf Proof of Lemma \ref{Lemma I4I7I5I8}}\\
\textit{\underline{For $\tilde I_4$.}} Let be consider the following random variable
\[
\tilde {\cal Z}_{m,k}(i) \equiv \e1_{{\cal U}_{m,k}(i) <m+k^0_m} \cro{ f(\X_{{\cal U}_{m,k}(i)}; \ebo)-f(\X_{{\cal U}_{m,k}(i)}; \hat \eb_{m+k})}-\frac{1}{m+k}\sum^{m+k^0_m}_{i=1} [f(\X_i;\ebo)-f(\X_i;\hat \eb_{m+k})] .
\]
It is obvious that $\eE^*_{m,k}[\tilde {\cal Z}_{m,k}(i)]=0$. Since $\tilde I_4=\sum^{m+l}_{i=m+1} \tilde {\cal Z}_{m,k}(i)$, we have that $\eE^*_{m,k}[\tilde I_4]=0$. For the conditional variance of $\tilde {\cal Z}_{m,k}(i)$  we have, using a quadratic Taylor expansion and the triangular inequality 
\[
Var^*_{m,k}[\tilde {\cal Z}_{m,k}(i)] = \eE^*\cro{\e1_{{\cal U}_{m,k}(i) <m+k^0_m}[ f(\X_{{\cal U}_{m,k}(i)}; \ebo)-f(\X_{{\cal U}_{m,k}(i)}; \hat \eb_{m+k})]}^2
=\frac{1}{m+k} \sum^{m+k^0_m}_{j=1} [f(\X_j;\ebo)-f(\X_j;\hat \eb_{m+k})]^2
\]
\[
=\frac{1}{m+k} \sum^{m+k^0_m}_{j=1} \cro{(\hat \eb_{m+k}-\ebo)^t\ef(\X_j;\ebo)+\frac{(\hat \eb_{m+k} -\ebo)^t}{2} \eff(\X_j; \tilde \eb_{m,k}) (\hat \eb_{m+k} -\ebo) }^2
\]
\[
\leq \frac{(m+k^0_m)\|\hat \eb_{m+k} -\ebo\|^2_2}{m+k}  \|\B_{m+k^0_m}\|^2_2+\frac{\|\hat \eb_{m+k} -\ebo\|^2_2}{m+k}  \left\| \sum^{m+k^0_m}_{j=1} \ef^t(\X_j;\ebo) \eff(\X_j; \tilde \eb_{m+k}) (\hat \eb_{m+k}-\ebo)\right\|_2 \]
\[\qquad \qquad +\frac{\|\hat \eb_{m+k} -\ebo\|^4_2}{4(m+k)} \sum^{m+k^0_m}_{j=1}\| \eff(\X_j; \tilde \eb_{m,k})\|^2_2.
\]
By the  Bienaym\'e-Tchebychev inequality we have  that
$(m+k)^{-1} \|\sum^{m+k^0_m}_{j=1} \ef^t(\X_j;\ebo) \eff(\X_j; \tilde \eb_{m+k}) \|_2 \leq (m+k)^{-1}  \sum^{m+k^0_m}_{j=1} \|\ef^t(\X_j;\ebo) \|_2 \| \eff(\X_j; \tilde \eb_{m+k})\|_2
\leq ( (m+k)^{-1}  \sum^{m+k^0_m}_{j=1} \ef^t(\X_j;\ebo) \ef(\X_j;\ebo))^{1/2}  ((m+k)^{-1}  \sum^{m+k^0_m}_{j=1} \| \eff(\X_j; \tilde \eb_{m+k})\|^2_2 )^{1/2}$.
Since $\Theta$ is compact, together the assumptions  (A2), (A4), we obtain that  $Var^*_{m,k}[\tilde {\cal Z}_{m,k}(i)] \leq C$. As in the linear case, using
$\sum^{T_m}_{l=1} g^{-2}(m,l,\gamma)  \leq C$,  we obtain the inequality (\ref{ineg_I4}) for $\tilde I_4$.\\
\textit{\underline{For $\tilde I_7$.}} We can prove in a similar way as in the  Lemma \ref{Lemma I3I6}, that $\eE^*_{m,k}[ g^{-1}(m,Thebychevl,\gamma) \tilde I_7 ]=o_{\eP}(1)$. Similar  arguments as for  $\tilde I_6$ of the Lemma \ref{Lemma I3I6}, unlike that $Var^*[\tilde I_7]= O_{\eP}(1)$ uniformly in $k$  and  probability 1, we obtain the relation (\ref{ineg_I4}) for $j=7$, by the  Bienaym\'e-Tchebychev inequality, Proposition \ref{Proposition A}, and inequality (\ref{CC}).\\
\textit{\underline{For $\tilde I_5$ and $\tilde I_8.$}} We consider the random variable  defined by 
\[
{\tilde{\tilde {\cal Z}}}_{m,k}(i) \equiv \e1_{{\cal U}_{m,k}(i) >m+k^0_m} \cro{ f(\X_{{\cal U}_{m,k}(i)}; \eb^0_m)-f(\X_{{\cal U}_{m,k}(i)}; \hat \eb_{m+k})}-\frac{1}{m+k} \sum^{m+k}_{j=m+k^0_m+1} [ f(\X_j;\eb^0_m)-f(\X_j;\hat \eb_{m+k})],
\]
and the results are proved by a similar way using assumption (A5) on the place of (A4).
\hspace*{\fill}$\blacksquare$\\

We recall the notations $
\hat \sigma^{(*)2}_{m,k}=(m-q)^{-1} \sum^{m}_{i=1} [\varepsilon^*_{m,k}(i)-m^{-1}\sum^m_{j=1} \ef^t(\X_j;\ebo)\varepsilon^*_{m,k}(j)) \B^{-1}_{m} \ef(\X_i;\ebo) ]^2
$ and  $\hat \sigma^2_{m,k}=(m+k)^{-1} \sum^{m+k}_{i=1}(\varepsilon_i-\bar \varepsilon_{m+k})^2$. We will prove that, under hypothesis $H_0$, $\hat \sigma^2_{m,k}$ and $\hat \sigma^{(*)2}_{m,k}$ are two uniformly consistent estimators for the variance $\sigma^2$ of the errors $\varepsilon^2$. Under hypothesis $H_1$, this two statistics are significantly different.
\begin{lemma}
\label{lemma 7}
 Suppose that the assumptions (A1)-(A4) hold.  \\
a) Under the  hypothesis $H_0$, we have, in probability,
\[
\sup_{1 \leq k \leq T_m} \eP^*_{m,k} \pth{ \left|\frac{\hat \sigma_{m,k}}{\hat \sigma^{(*)}_{m,k}} -1\right| \geq \epsilon } {\underset{m \rightarrow \infty}{\longrightarrow}} 0.
\]
b) If furthermore the assumption (A5) holds, under the  hypothesis $H_1$, for all $\epsilon >0$, there exists a constant $M>0$ such that
\[
\sup_{1 \leq k \leq T_m} \eP^*_{m,k} \pth{ \left|\frac{\hat \sigma_{m,k}}{\hat \sigma^{(*)}_{m,k}} -1\right| \geq M } \leq \epsilon+o_{\eP}(1) .
\]
\end{lemma}
\noindent {\bf Proof of Lemma \ref{lemma 7}}\\
We have the decomposition, for each  $i$ of 1 to $m$: $\varepsilon^*_{m,k}(i)-(m^{-1}\sum^m_{j=1} \ef^t(\X_j;\ebo)\varepsilon^*_{m,k}(j)) \B^{-1}_{m} \ef(\X_i;\ebo) \equiv \sum^8_{j=1}J_j(m,k,i)$,
where: $J_1(m,k,i)=\varepsilon_{{\cal U}_{m,k}(i)}$, $J_2(m,k,i)=-\pth{m^{-1}\sum^m_{j=1} \ef^t(\X_j;\ebo)\varepsilon_{{\cal U}_{m,k}(j)}} \B^{-1}_{m} \ef(\X_i;\ebo)$,
\[
J_3(m,k,i)=\e1_{k \leq k^0_m} \e1_{{\cal U}_{m,k}(i) < m+k^0_m}\left[f(\X_{{\cal U}_{m,k}(i)};\ebo)-f(\X_{{\cal U}_{m,k}(i)};\hat \eb_{m+k})   -\frac{1}{m+k} \sum^{m+k}_{j=1}[f(\X_j;\ebo)-f(\X_j;\hat \eb_{m+k})]\right],
\]
\[
J_4(m,k,i)= \e1_{k > k^0_m} \e1_{{\cal U}_{m,k}(i) < m+k^0_m}\cro{f(\X_{{\cal U}_{m,k}(i)};\ebo)-f(\X_{{\cal U}_{m,k}(i)};\hat \eb_{m+k}) }  -\e1_{k >k^0_m} \frac{1}{m+k} \sum^{m+k^0_m}_{j=1} [ f(\X_j;\ebo)-f(\X_j;\hat \eb_{m+k})],
\]
\[
J_5(m,k,i)=\e1_{k > k^0_m} \e1_{{\cal U}_{m,k}(i) > m+k^0_m}\cro{f(\X_{{\cal U}_{m,k}(i)};\eb^0_m)-f((\X_{{\cal U}_{m,k}(i)};\hat \eb_{m+k})) }-\e1_{k >k^0_m} \frac{1}{m+k} \sum^{m+k}_{j=m+k^0_m+1} [ f(\X_j;\eb^0_m)-f(\X_j;\hat \eb_{m+k})],
\]
\[
J_6(m,k,i)=- \pth{\frac{1}{m} \sum^m_{j=1} \ef^t(\X_j;\ebo)\e1_{k \leq k^0_m} \e1_{{\cal U}_{m,k}(j) < m+k^0_m}\cro{f(\X_{{\cal U}_{m,k}(j)};\ebo)-f((\X_{{\cal U}_{m,k}(j)};\hat \eb_{m+k})) }}\B^{-1}_{m}\ef(\X_i;\ebo) \]
\[\qquad \qquad \qquad \qquad  +\frac{1}{m+k} \sum^{m+k}_{j=1}[f(\X_j;\ebo)-f(\X_j;\hat \eb_{m+k})],
\]
\[
J_7(m,k,i)=- \pth{\frac{1}{m} \sum^m_{j=1} \ef^t(\X_j;\ebo)\e1_{k > k^0_m} \e1_{{\cal U}_{m,k}(j) < m+k^0_m}\cro{f(\X_{{\cal U}_{m,k}(j)};\ebo)-f(\X_{{\cal U}_{m,k}(j)};\hat \eb_{m+k}) }}  \B^{-1}_{m}\ef(\X_i;\ebo)
\]
\[
\qquad \qquad \qquad \qquad   +\e1_{k >k^0_m} \frac{1}{m+k} \sum^{m+k^0_m}_{j=1} [ f(\X_j;\ebo)-f(\X_j;\hat \eb_{m+k})],
\]
\[
J_8(m,k,i)=- \pth{\frac{1}{m} \sum^m_{j=1} \ef^t(\X_j;\ebo)\e1_{k > k^0_m} \e1_{{\cal U}_{m,k}(j) > m+k^0_m}\cro{f(\X_{{\cal U}_{m,k}(j)};\eb^0_m)-f(\X_{{\cal U}_{m,k}(j)};\hat \eb_{m+k}) }} \B^{-1}_{m}\ef(\X_i;\ebo)
\]
\[
\qquad \qquad \qquad \qquad  +\e1_{k >k^0_m} \frac{1}{m+k} \sum^{m+k}_{j=m+k^0_m+1} [ f(\X_j;\eb^0_m)-f(\X_j;\hat \eb_{m+k})].
\]
Then 
\begin{equation}
\label{smk}
\hat \sigma^{(*)2}_{m,k}=\frac{1}{m-q} \sum^m_{i=1} \cro{ \sum^8_{j=1}J_j(m,k,i)}^2.
\end{equation}
Under $H_0$, we have $J_4,J_5, J_7,J_8=0$.\\
Following results hold under the two hypotheses $H_0$ and $H_1$. For $J_1$ we have that for any $\epsilon >0$
\begin{equation}
\label{J1}
\sup_{1 \leq k} \eP^*_{m,k} \cro{ \left|\frac{1}{m-q} \frac{1}{\hat \sigma^2_{m,k}} \sum^m_{i=1} J^2_1(m,k,i) -1 \right| \geq \epsilon} {\underset{m \rightarrow \infty}{\longrightarrow}} 0.
\end{equation}
 For $J_2$, we have $\eE^*_{m,k}[ J_2(m,k,i)]=- \bar \varepsilon_{m+k} m^{-1}\sum^m_{j=1} \ef^t(\X_j;\ebo) \B^{-1}_{m} \ef(\X_i;\ebo)$. 
Then, using the independence of $\varepsilon_i$ and of $\X_i$, assumption (A4), we obtain the  convergence in probability, uniformly in $k$, as $m \rightarrow \infty$,
\[
\eE^*_{m,k}[ \frac{1}{m} \sum^m_{i=1}J_2(m,k,i)]=- \bar \varepsilon_{m+k} \pth{\frac{1}{m}\sum^m_{j=1} \ef^t(\X_j;\ebo) \B^{-1}_{m} \frac{1}{m} \sum^m_{i=1} \ef(\X_i;\ebo)} \rightarrow 0 \cdot \A^t \B^{-1} \A= 0.
\]
We have also the approximation of $\eE^*_{m,k}[ m^{-1} \sum^m_{i=1}J_2^2(m,k,i)]$ by
$ \B^{-1}_{m }\eE^*_{m,k} [ m^{-2} \sum^m_{j=1} \ef^t(\X_j;\ebo)$ $\cdot \ef(\X_j;\ebo)\varepsilon^2_{{\cal U}_{m,k}(j)}$\\$+2m^{-2} \sum_j \sum_{j' \neq j} \ef^t(\X_j;\ebo) \ef(\X_{j'};\ebo) \varepsilon_{{\cal U}_{m,k}(j)} \varepsilon_{{\cal U}_{m,k}(j')}]$
$=  \B^{-1}_{m } [ m^{-1} \B_m (m+k)^{-1} \sum^{m+k}_{j=1} \varepsilon^2_j+2m^{-2} \sum_j \sum_{j' \neq j} \ef(\X_j;\ebo)$\\$ \cdot \ef^t(\X_{j'};\ebo)  ((m+k)^{-1} \sum^{m+k}_{a=1} \varepsilon_a)^2 ]$.
But $(m+k)^{-1} \sum^{m+k}_{a=1} \varepsilon_a  \overset{{\eP}} {\underset{m \rightarrow \infty}{\longrightarrow}} 0$ and $(m+k)^{-1} \sum^{m+k}_{a=1} \varepsilon^2_a  \overset{{\eP}} {\underset{m \rightarrow \infty}{\longrightarrow}} \sigma^2$. Hence,  we have uniformly in $k$
\begin{equation}
\label{b1}
\eE^*_{m,k}[ m^{-1} \sum^m_{i=1}J_2^2(m,k,i)]$ $ \overset{{\eP}} {\underset{m \rightarrow \infty}{\longrightarrow}} 0.
\end{equation}
On the other hand, we can write  $J_3(m,k,i)={\cal Z}_{k,m}(i)$, with ${\cal Z}_{k,m}(i)$ defined by the relation (\ref{zz}). By the proof of the Lemma \ref{Lemma I3I6}, since $Var^*_{m,k}[{\cal Z}_{k,m}(i)] = (m+k)^{-1}\B(1+o_{\eP}(1))$ and the Bienaym\'e-Tchebychev inequality, we have in  probability
\begin{equation}
\label{b2}
\sup_{k \geq 1} \eP^*_{m,k} \cro{\frac{1}{m-q} \sum^m_{i=1} J_3^2(m,k,i) \geq \epsilon} {\underset{m \rightarrow \infty}{\longrightarrow}} 0.
\end{equation}
For $J_1J_2$ we have
$
m^{-1} \sum^m_{i=1} J_1(m,k,i) J_2(m,k,i)=- \pth{m^{-1} \sum^m_{j=1} \ef^t(\X_j;\ebo) \varepsilon_{{\cal U}_{m,k}(j)}} \B^{-1}_{m}\pth{m^{-1} \sum^m_{j=1} \ef(\X_j;\ebo) \varepsilon_{{\cal U}_{m,k}(j)}}$
and then, as for the calculations from above for $J_2$, we have $\eE^{*}_{m,k}[m^{-1} \sum^m_{i=1} J_1(m,k,i) J_2(m,k,i) ]=o_{\eP}(1)$, uniformly in $k$.\\
For $J_1J_3$ we have, using the assumptions (A2)-(A4), 
$
\eE^{*}_{m,k}[m^{-1}  \sum^m_{i=1} J_1(m,k,i) J_3(m,k,i) ]=m^{-1}  \sum^m_{i=1} \{(m+k)^{-1} \sum^{m+k}_{j=1} \varepsilon_j [f(\X_j;\ebo)-f(\X_j;\hat \eb_{m+k})]$.
$
- \bar \varepsilon_{m+k}(m+k)^{-1} \sum^{m+k}_{j=1} \varepsilon_j [f(\X_j;\ebo)-f(\X_j;\hat \eb_{m+k})]\} =o_{\eP}(1)-o_{\eP}(1)o_{\eP}(1)=o_{\eP}(1)$. 
We show similar for the other cases that $\eE^{*}_{m,k}[m^{-1} \sum^m_{i=1} J_1(m,k,i) J_l(m,k,i) ]\overset{{\eP}} {\underset{m \rightarrow \infty}{\longrightarrow}} 0$, $l=3,4, \cdots, 8$. 
The conditional  expectation  $\eE^{*}_{m,k}[m^{-1} \sum^m_{i=1} J_2(m,k,i) J_3(m,k,i) ]$ is equal to
\[
- \frac{1}{m} \sum^m_{i=1} \ef^t(\X_i;\ebo) \frac{1}{m}
\frac{1}{m+k} \sum^{m+k}_{j=1} \varepsilon_j [f(\X_j;\ebo)-f(\X_j;\hat \eb_{m+k}) ] \B^{-1}_{m} \ef(\X_j;\ebo)\]
\[\qquad \qquad \qquad +\bar \varepsilon_{m+k} \pth{\frac{1}{m} \sum^m_{j=1} \ef^t(\X_j;\ebo)  } \B^{-1}_{m} \pth{ \frac{1}{m} \sum^m_{i=1} \ef(\X_i;\ebo) }\frac{1}{m+k} \sum^{m+k}_{j=1}[f(\X_j;\ebo)-f(\X_j;\hat \eb_{m+k}) ], 
\]
which  converges to 0, from $m \rightarrow \infty$, uniformly in $k$, since $\varepsilon_i$ is independent of $\X_i$, $\bar \varepsilon_{m+k} \rightarrow 0$, and all terms in $\X_i$ are bounded.\\
\textit{a) } \underline{Under hypothesis $H_0$}. For $J_6(m,k,i)$ we have
\[  
\eE^*_{m,k}[ J_6(m,k,i)]=- (\hat \eb_{m+k} -\ebo)^t \pth{\frac{1}{m+k} \sum^{m+k}_{a=1} \ef(\X_a;\ebo)} \pth{\frac{1}{m} \sum^m_{j=1} \ef^t(\X_j;\ebo)}\B^{-1}_{m}\ef(\X_i;\ebo).
\]
Then, under assumption (A4),
$\eE^*_{m,k}[  m^{-1} \sum^m_{i=1}J_6(m,k,i)]=-(\hat \eb_{m+k} -\ebo)^t \A \cdot \A^t \cdot \B^{-1} \A (1+o_{\eP}(1))$. 
Similarly, with $C$ a  constant, $\eE^{*}_{m,k}[ J^2_6(m,k,i)]\leq \|\hat \eb_{m+k} -\ebo\|_2^2\cdot  \|\A \A^t \|_2^2 \cdot \|\B^{-1}\|^2_2 \cdot  \ef(\X_i;\ebo) \ef^t(\X_i;\ebo) (1+o_{\eP}(1))$, thus
$\eE^*_{m,k}[ m^{-1}\sum^m_{i=1}J^2_6(m,k,i)]\leq \|\hat \eb_{m+k} -\ebo\|_2^2 \cdot  \|\A \A^t \|_2^2 \cdot \|\B^{-1}\|^2_2 \cdot  \|\B\|_2  \cdot (1+o_{\eP}(1))$ and 
$\eE^*_{m,k}[ m^{-1}\sum^m_{i=1}J^4_6(m,k,i)]=C \|\hat \eb_{m+k} -\ebo\|_2^4  (1+o_{\eP}(1))$. 
Hence, under $H_0$,  $Var^*_{m,k}[ m^{-1}\sum^m_{i=1}J^2_6(m,k,i)]=o_{\eP}(1)$. By the Bienaym\'e-Tchebychev inequality, we have in  probability:
\[
\sup_{k \geq 1} \eP^*_{m,k} \cro{\frac{1}{m-q} \sum^m_{i=1} J_6^2(m,k,i) \geq \epsilon} {\underset{m \rightarrow \infty}{\longrightarrow}} 0. 
\]
The relations (\ref{smk}), (\ref{J1}) and since all other conditional expectations for   $\hat \sigma^{(*)2}_{m,k}$ expression are negligible, imply the assertion \textit{(a)}. \\
\textit{b) }\underline{Under $H_1$.}
For $b \in \{4,5,6,7,8 \}$ we will prove that for all  $\epsilon >0$ there exists a  $M>0$ such that
\begin{equation}
\label{J4}
\sup_{k\geq 1} \eP^*_{m,k} \cro{\frac{1}{m-q}\sum^m_{i=1}J^2_b(m,k,i)\geq M} \leq \epsilon+ o_{\eP}(1).
\end{equation}
 In view of  the previous calculus for $J_6$ we have that the relation (\ref{J4}) holds for $b=6$. \\
  For $J_4(m,k,i)$ we have $J_4(m,k,i)=\tilde {\cal Z}_{k,m}(i)$ and by the proof of the Lemma \ref{Lemma I4I7I5I8} and the  Bienaym\'e-Tchebychev inequality, we have that  the relation (\ref{J4}) holds for $b=4$. \\
For $J_7$, its conditional expectation  $\eE^{*}_{m,k}[ J_7(m,k,i)]$ is, by Taylor expansions and using assumptions (A2), (A4), 
\[
 \frac{(\ebo- \hat \eb_{m+k})^t}{m+k} \sum^{m+k^0_m}_{a=1} \cro{\ef(\X_a;\ebo)+\eff(\X_a;\tilde \eb) \frac{\ebo- \hat \eb_{m+k}}{2}  } \cro{ 1- \frac{1}{m} B^{-1}_{m} \sum^m_{j=1} \ef(\X_j;\ebo)} \ef^t(\X_i;\ebo)(1+ o_{\eP}(1)).
\]
Similarly $\eE^{*}_{m,k}[  m^{-1}\sum^m_{i=1} J_7^2(m,k,i)]=C \|\ebo- \hat \eb_{m+k}\|_2^2(1+o_{\eP}(1))$, $\eE^{*}_{m,k}[  m^{-1}\sum^m_{i=1} J_7^4(m,k,i)]=C \|\ebo- \hat \eb_{m+k}\|_2^4(1+o_{\eP}(1))$ 
which imply the relation  (\ref{J4}) for $J_7$. By Bienaym\'e-Tchebychev inequality, we obtain  (\ref{J4}) for $J_5$ and $J_8$ using assumption (A5) on the place of (A4). \\
Now we consider the product of the terms of different suffix.
The products of  $J_3$ with $J_4, J_5, J_7,J_8$, of $J_4$ with $J_5, J_6, J_8$ and of $J_5$ with $J_6,J_7$ are 0.
For $J_1J_4$ we have $\eE^{*}_{m,k}[m^{-1} \sum^m_{i=1} J_1(m,k,i) J_4(m,k,i) ]=(m+k)^{-1} \sum^{m+k^0_m}_{j=1}[\varepsilon_j-\bar \varepsilon_{m+k}] [f(\X_j;\ebo)-f(\X_j;\hat \eb_{m+k})]$. But, for all $\epsilon>0$ there exists $M>0$ such that $
\eP\cro{(m+k)^{-1} \sum^{m+k^0_m}_{j=1} |f(\X_j;\ebo)-f(\X_j;\hat \eb_{m+k}) | >M} <\epsilon$, from which, together the fact $\eE[\varepsilon_j]=0$,  one may deduce that\\  $\eE^{*}_{m,k}[m^{-1} \sum^m_{i=1} J_1(m,k,i) J_4(m,k,i) ]\overset{{\eP}} {\underset{m \rightarrow \infty}{\longrightarrow}} 0 $, uniformly in $k$. 
By similar arguments we prove the uniformly convergence to 0 in probability,  for  all other combinations of $J_2$ and $J_4, \cdots, J_8$. The  not insignificant terms are $J^2_4, J^2_7, J^2_5, J^2_8$. We consider now $J_4J_7$, the other cases are similar. Taking into account the fact that  $\eE^*_{m,k}[f(\X_{{\cal U}_{m,k}(i)};\eb)f(\X_{{\cal U}_{m,k}(j)};\eb)]$ is  equal to  $\eE^*_{m,k}[f(\X_{{\cal U}_{m,k}(i)};\eb)]\eE^*[f(\X_{{\cal U}_{m,k}(j)};\eb)]$ for $i \neq j$ and to  $\eE^*[f^2(\X_{{\cal U}_{m,k}(i)};\eb)]$ for $i=j$, we have that $\eE^{*}_{m,k}[m^{-1} \sum^m_{i=1} J_4(m,k,i) J_7(m,k,i) ]=$
\[
\frac{1}{m(m+k)}\sum^{m+k^0_m}_{j=1}[f(\X_j;\ebo)-f(\X_j;\hat \eb_{m+k}) ]^2 \pth{\frac{1}{m} \sum^m_{i=1}\ef^t(\X_i;\ebo)}\B^{-1}_{m} \pth{\frac{1}{m} \sum^m_{i=1}\ef(\X_i;\ebo)}(1+o_{\eP}(1)),
\]
which converges to 0 in probability, uniformly in $k$.\\
Hence, in conclusion, taking into account the relations (\ref{smk}), (\ref{b1}), (\ref{b2}) and (\ref{J4}), 
\[
\hat \sigma_{m,k}^{2(*)}=\frac{1}{m-q} \cro{\sum^m_{i=1}(J^2_1(m,k,i)+J^2_4(m,k,i)+J^2_5(m,k,i)+J^2_7(m,k,i)+J^2_8(m,k,i))} (1+o_{\eP}(1))\]
\[
\geq \frac{1}{m}\sum^m_{i=1}J^2_1(m,k,i)(1+o_{\eP}(1))
\]
and the assertion \textit{(b)} follows by (\ref{J1}). 
\hspace*{\fill}$\blacksquare$\\

\end{document}